\documentclass[10pt,a4paper]{amsart}

\usepackage{amsmath,amsfonts,amssymb,amsthm,latexsym,amscd}
\usepackage[all]{xy}
\usepackage[mathscr]{eucal}

\theoremstyle{plain}
\newtheorem{thm}{Theorem}[section]
\newtheorem{prop}[thm]{Proposition}

\newtheorem{lem}[thm]{Lemma}

\theoremstyle{definition}
\newtheorem{defn}[thm]{Definition}

\theoremstyle{remark}
\newtheorem{rem}[thm]{Remark}
\newtheorem{expl}[thm]{Example}

\numberwithin{equation}{section}

\newcommand{\sD}{\mathcal{D}}

\newcommand{\sJ}{\mathcal{J}}

\newcommand{\sL}{\mathcal{L}}
\newcommand{\sN}{\mathcal{N}}
\newcommand{\sS}{\mathcal{S}}
\newcommand{\sP}{\mathcal{P}}
\newcommand{\sO}{\mathcal{O}}

\newcommand{\sStw}{\widetilde{\sS}}

\newcommand{\bL}{\mathbf{L}}
\newcommand{\bS}{\mathbf{S}}

\renewcommand{\AA}{\mathbb{A}}
\newcommand{\BB}{\mathbb{B}}

\newcommand{\DD}{\mathbb{D}}

\newcommand{\RR}{\mathbb{R}}
\renewcommand{\SS}{\mathbb{S}}
\newcommand{\ZZ}{\mathbb{Z}}

\newcommand{\id}{\textup{id}}
\newcommand{\im}{\textup{im}}
\newcommand{\Hom}{\textup{Hom}}

\newcommand{\nin}{\noindent}

\newcommand{\Spaces}{\textup{Spaces}}

\newcommand{\ra}{\rightarrow}
\newcommand{\lra}{\longrightarrow}
\newcommand{\co}{\colon\!}

\newcommand{\iso}{^{\textup{iso}}}
\newcommand{\uli}{\underline}

\newcommand{\hocolim}{\textup{hocolim}}
\newcommand{\colim}{\textup{colim}}
\newcommand{\mor}{\textup{mor}}
\newcommand{\ind}{\textup{ind}}

\newcommand{\op}{^{\textup{op}}}
\newcommand{\cat}{\textup{cat}}

\newcommand{\ddiamond}{{\diamond\diamond}}

\newcommand{\Or}{\textup{O}}
\newcommand{\hofiber}{\textup{hofiber}}

\newcommand{\smin}{\smallsetminus}

\newcommand{\ee}{\'{e} }

\newcommand{\colimsub}[1]{\begin{array}[t]{cc} \textup{colim} \\
[-1.7mm] \scriptstyle{#1} \end{array}}

\newcommand{\hocolimsub}[1]{\begin{array}[t]{cc} \textup{hocolim} \\
[-1mm] \scriptstyle{#1} \end{array}}

\newcommand{\twosub}[2]{\begin{array}{cc}
\scriptstyle{#1} \\  [-1mm] \scriptstyle{#2}  \end{array}}

\begin{document}

\title[The Block Structure Spaces and Orthogonal Calculus]{The Block Structure
Spaces of Real Projective Spaces and Orthogonal Calculus of Functors
II}

\author{Tibor Macko and Michael Weiss}

\date{\today}



\address{Mathematisches Institut \\ Universit\"at M\"unster \\
Einsteinstrasse 62 \\ M\"unster, D-48149 \\ Germany \\ and
Matematick\'y \'Ustav SAV \\ \v Stef\'anikova 49 \\ Bratislava,
SK-81473 \\ Slovakia} \email{macko@math.uni-muenster.de}
\address{Department of Mathematical Sciences \\ University of Aberdeen
\\ Aberdeen, AB24 3UE \\ Scotland, UK}
\email{mweiss@maths.abdn.ac.uk}

\thanks{M Weiss supported by the Royal Society (Wolfson
research merit award)}

\begin{abstract} For a finite dimensional real vector space $V$ with
inner product, let $F(V)$ be the block structure space, in the sense
of surgery theory, of the projective space of $V$. Continuing a
program launched in \cite{Ma}, we investigate $F$ as a functor on
vector spaces with inner product, relying on functor calculus ideas.
It was shown in \cite{Ma} that $F$ agrees with its first Taylor
approximation $T_1F$ (which is a polynomial functor of degree 1) on
vector spaces $V$ with $\dim(V)\ge 6$. To convert this theorem into
a functorial homotopy-theoretic description of $F(V)$, one needs to
know in addition what $T_1F(V)$ is when $V=0$. Here we show that
$T_1F(0)$ is the standard $L$-theory space associated with the group
$\ZZ/2$, except for a deviation in $\pi_0$. The main corollary is a
functorial two-stage decomposition of $F(V)$ for $\dim(V)\ge 6$
which has the $L$-theory of the group $\ZZ/2$ as one layer, and a
form of unreduced homology of $\RR P(V)$ with coefficients in the
$L$-theory of the trivial group as the other layer. (Except for
dimension shifts, these are also the layers in the traditional
Sullivan-Wall-Quinn-Ranicki decomposition of $F(V)$. But the
dimension shifts are serious and the SWQR decomposition of $F(V)$ is
not functorial in $V$.) Because of the functoriality, our analysis
of $F(V)$ remains meaningful and valid when $V=\RR^\infty$.
\end{abstract}

\maketitle

\section{Introduction}
This paper is a continuation of \cite{Ma}. In \cite{Ma} a certain
continuous functor from the category $\sJ$ of finite-dimensional
real vector spaces with inner product to the category $\Spaces_\ast$
of pointed spaces was introduced. Here (in this introduction) we
denote this functor by $F^g$. For $V$ an object of $\sJ$, the value
$F^g(V)$ is the block structure space
\[ \sStw (\RR P(V)) \]
of the projective space of $V$. For a morphism $\xi$ in $\sJ$ the
map $F^g(\xi)$ is a generalization of the join construction of Wall.
See \cite{Ma} or \cite{Qu} for the definition of the block structure
space of a manifold, \cite{Ma} and \cite{Wa} for more on the join
construction.

Each space $F^g(V)$, alias block structure space of $\RR P(V)$, is
individually well understood as the $n$-fold loop space of the
homotopy fiber of a standard assembly map in $L$-theory, where
$n=\dim(V)-1$ and we assume $n\ge 5$. See \cite{Qu}. The assembly
map has good naturality properties, but the prefix $\Omega^n$ tends
to corrupt these when $n$ becomes a variable. Hence the standard
methods for calculating the values $F^g(V)$ do not lead to very
satisfying homotopy theoretic descriptions of the induced maps
$F^g(\xi)$.

The goal of the project presented in \cite{Ma} and here is to
provide a homotopy theoretic description of the spaces $F^g(V)$
natural in $V$, i.e. to describe the functor $F^g$. This will allow
us to let $\dim (V)$ tend to infinity. Hence it gives us a homotopy
theoretic description of the (homotopy) colimit of the spaces
$F^g(\RR^n)$, a space which might be considered as the block
structure space of $\RR P^\infty$ and which, as explained in the
introduction to \cite{Ma}, can be further related to certain space
of equivariant (honest = non-block) automorphisms of spheres. For
this purpose another tool is employed, the orthogonal calculus of
functors of Weiss \cite{We}. The desired description should be
obtained from the orthogonal calculus Taylor tower of the functor
$F^g$. This tower yields in particular a first Taylor approximation
$T_1F^g$ of $F^g$ ``at infinity'', which is another functor from
$\sJ$ to $\Spaces_\ast$ and comes with a canonical transformation
$F^g\to T_1F^g$. The degree 1 property of $T_1F^g$ implies a
homotopy fibration sequence, natural in $V$:
\begin{equation}\label{alt-des}
\Omega^{\infty}[ (S(V)_+ \wedge \mathbf{\Theta}^{(1)})_{h\Or(1)} ]
\rightarrow T_1F^g(0) \rightarrow T_1F^g(V).
\end{equation}
Here $S(V)$ is the unit sphere in $V$ with the antipodal involution,
the subscript $+$ denotes an added base point,
$\mathbf{\Theta}^{(1)}$ denotes the first derivative spectrum of
$F^g$ and the subscript $h\Or(1)$ denotes a homotopy orbit
construction for the symmetry group $\Or(1)\cong\ZZ_2$. See
\cite{We} or \cite{Ma} for the definitions and more on the Taylor
tower of a continuous functor from $\sJ$ to $\Spaces_\ast$.

A first step in the project was made in \cite{Ma}. Namely, it was
shown that for $V$ such that $\dim (V) \geq 6$ the canonical map
$F^g(V) \rightarrow T_1F^g(V)$ is a homotopy equivalence. Therefore
the homotopy fibration sequence (\ref{alt-des}) can be rewritten as
a homotopy fibration sequence
\begin{equation}\label{alt-des-2}
\Omega^{\infty}[ (S(V)_+ \wedge \mathbf{\Theta}^{(1)})_{h\Or(1)} ]
\rightarrow T_1F^g(0) \rightarrow F^g(V).
\end{equation}
It was also shown in \cite{Ma} that $\pi_k\mathbf{\Theta}^{(1)}$ is
the $k$-th $L$-group of the trivial group. We see that in order to
obtain the desired description of $F^g(V)$ when $\dim(V)\ge 6$,
natural in $V$, it is (nearly) enough to describe $T_1F^g(0)$.
Unfortunately this is far from easy.

Let $\AA$ be the additive category of finitely generated free
abelian groups and let $\AA[\ZZ_2]$ be the additive category of
finitely generated free modules over the group ring $\ZZ[\ZZ_2]$.
Both of these have duality functors and hence determine $L$-groups
and $L$-theory spaces. Our main result is as follows.

\vspace{2mm} \noindent \textbf{Theorem 2.} {\it There is a homotopy
fibration sequence
\[ \CD T_1 F^g(0) @>>> \bL_0(\AA[\ZZ_2]^+)
@>{\tilde\sigma/8}>> \ZZ\,. \endCD\]
} \vspace{2mm} 

The space $\bL_0(\AA[\ZZ_2]^+)$ is the standard $L$-theory space for
the group $\ZZ_2$ with the trivial orientation character. Its
homotopy groups are $\pi_k\bL_0(\AA[\ZZ_2]^+) = L_k(\ZZ_2^+)$. The
map $\tilde\sigma/8$ is the composition of the transfer homomorphism
$L_0\AA[\ZZ_2]^+)\to L_0(\AA)$ and the isomorphism
$L_0(\AA)\cong\ZZ$ defined by signature over $8$. It is onto. See
\cite[chapter 13A]{Wa}.

\medskip
In \cite{Ma} mainly the geometric surgery theory of Wall \cite{Wa}
was used, but for the proof of theorem 2 we switch to the algebraic
theory of surgery of Ranicki \cite{Ra}. In the geometric surgery
setup, the block structure space of an $n$-dimensional closed
manifold $X$ fits into the homotopy fibration sequence, due to Quinn
\cite{Qu},
\begin{equation} \label{sur-seq-intr}
\sStw (X) \rightarrow \sN (X) \rightarrow \sL_n (X),
\end{equation}
where $\sN (X)$ is the space of normal invariants of $X$ and $\sL_n
(X)$ is the surgery obstruction space associated with $X$. The
homotopy groups of $\sL_n (X)$ are the $L$-groups of $\pi = \pi_1
(X)$ with associated orientation character and with an appropriate
dimension shift. Quinn's homotopy fibration sequence is the space
version of a long exact sequence of (homotopy) groups usually
attributed to Sullivan in the simply connected case, and to Wall in
the nonsimply connected case.

In the algebraic surgery setup the input consists typically of a
connected simplicial complex $X$, an integer $n$, a universal
covering space for $X$ with deck transformation group $\pi$ and a
homomorphism $w\co \pi\to \ZZ_2$. The output is the homotopy
fibration sequence, due to Ranicki \cite{Ra},
\begin{equation} \label{alg-sur-seq-intr}
\bS (X,n,w) \rightarrow \bL_n (\AA_\ast (X,w) ) \rightarrow \bL_n
(\AA[\pi]^w).
\end{equation}
Here all three spaces are $L$-theoretic spaces constructed from
certain additive categories with chain duality. The homotopy groups
$\pi_k\bL_n(\AA[\pi]^w)$ are again the groups $L_{k+n}$ of the group
ring $\ZZ[\pi]$ with the $w$-twisted involution. The homotopy groups
$\pi_k\bL_n (\AA_\ast (X,w) )$ are better known as the generalized
homology groups $H_{k+n}$ of $X$ with $w$-twisted coefficients in
the $L$-theory spectrum of the trivial group. We should perhaps add
that our conventions here are slightly less restrictive than the
ones which Ranicki has in many of his papers on algebraic surgery.
As a consequence we have strict $4$-periodicity, $\bS(X,n,w)\simeq
\Omega^4\bS(X,n,w)$ for all $n\in \ZZ$, in addition to the spectrum
property $\bS(X,n+1,w)\simeq \Omega\bS(X,n,w)$ which Ranicki also
insists on in most circumstances.

If $X$ is a triangulated closed $n$-manifold and $w$ is its
orientation character, then modulo a small modification the sequence
(\ref{alg-sur-seq-intr}) can be identified with the sequence
(\ref{sur-seq-intr}) by a result of \cite{Ra}. (In this situation we
usually write $\bS (X)$ instead of $\bS(X,n,w)$.) More concretely
the first terms in the sequences (the two versions of the block
structure space) are related via a homotopy fibration sequence
\begin{equation} \label{geo-alg-sur-comp}
\sStw (X) \ra \bS(X) \ra \ZZ.
\end{equation}   
The advantage of the algebraic setup is that $\bS(X)$ is much more
tractable from the point of view of algebraic topology: it is an
infinite loop space and it is 4-periodic as such, almost by
definition.

But for us the algebraic setup has some disadvantages, too. No
triangulation is invariant under the action of the orthogonal group
$\Or(V)$ on $V$. Therefore it is not possible to define a continuous
functor from $\sJ$ to $\Spaces_\ast$ by a formula such as $V \mapsto
\bS (\RR P(V))$. Instead we construct a continuous functor $F^a$
from $\sJ$ to $\Spaces_\ast$ along the following lines. For $V$ in
$\sJ$, the value $F^a(V)$ is a colimit of spaces $\bS(X)$, where $X$
runs through a directed system of generalized simplicial complexes
obtained from certain generalized triangulations of $\RR P(V)$. Of
course, the space $F^a(V)$ will have the homotopy type of $\bS (\RR
P(V))$.

We emphasize that, although $F^a$ is better behaved than $F^g$ from
the orthogonal calculus point of view, its behavior on objects $V$
in $\sJ$ of dimension $<3$ is still not good. This is due to the
fact that the map $\RR P(V)\to \RR P(W)$ induced by a morphism $V\to
W$ in $\sJ$ need not be $1$-connected if $\dim(V)<3$. However, that
difficulty can be overcome and we have the following result which
easily implies Theorem 2:

\vspace{2mm} \noindent \textbf{Theorem 1.} {\it We have $T_1 F^a(0)
\simeq \bL_0(\AA[\ZZ_2]^+)$. } \vspace{1mm}

An overview of the proof will be given in the next section. --- The
algebraic surgery approach also gives us the following:

\vspace{1mm} \noindent \textbf{Remark.} {\it The sequence {\rm
(\ref{alt-des-2})} is a homotopy fibration sequence of infinite loop
spaces. } \vspace{1mm}

The paper is organized as follows. In section \ref{proof} we give
the proof of our main result, Theorem 1, modulo certain theorems A,
B and C and a natural homotopy fiber sequence $F^g(V)\to F^a(V)\to
\ZZ$ when $\dim(V)\ge 6$. In the rest of the paper we prove theorems
$A$, $B$ and $C$ and construct that natural homotopy fiber sequence.
Specifically, in section~\ref{alg-sur} we give a review of the tools
from algebraic surgery we need. Section~\ref{oc-and-prod} contains a
somewhat abstract preview of the functor $F^a$ while
sections~\ref{alg-join} and~\ref{str-spaces} deliver the technical
details. Section~\ref{proofs} contains the proofs of Theorems A, B
and C. In section~\ref{geo-alg-join} we relate $F^a$ to $F^g$. This
completes the proof of theorem 1 and it also reduces Theorem 2 to
Theorem 1. At the very end of section~\ref{geo-alg-join} we also
explain the remark above on infinite loop space structures.


\section{Proof of Theorem 1.} \label{proof}

The main ingredients in the proof of Theorem 1 are Theorems A, B and
C below. Their statements and proofs rely on the algebraic theory of
surgery setup. The construction of the functor $F^a$ and the proofs
of Theorems A and B are themes of the subsequent sections. Theorem C
was proved in \cite{Ma}. Although $F^a$ has a complicated
definition, for the purposes of this section we may pretend that it
is given by $F^a \co V \mapsto \bS (\RR P(V))$. Now we state
Theorems A and B without proofs and derive Theorem 1 from them.

\vspace{2mm} \noindent \textbf{Theorem A.} {\it For an oriented $W
\in \sJ$ such that $4$ divides $\dim (W)$, there is a homotopy
equivalence $F^a(V) \ra \Omega^WF^a(V)$, natural in $V$. }
\vspace{1mm}

This is just the usual $4$-periodicity in $L$-theory. --- Our next
statement is about a Thom isomorphism in algebraic $L$-theory. Let
$V$ and $W$ objects in $\sJ$. As one would expect, there is a join
map $\bS(\RR P(W)) \ra \bS(\RR P(V\oplus W))$ and we do not
generally have a way of extending that to a map between the two
homotopy fiber sequences (\ref{alg-sur-seq-intr}) for $X=\RR P(W)$
and $X=\RR P(V\oplus W)$, respectively.  But it is relatively easy
to supply the lower horizontal arrow in a homotopy commutative
diagram
\[
\xymatrix{
  \bS(\RR P(W)) \ar[r] \ar[d] & \bS(\RR P(V\oplus W)) \ar[d] \\
\bL_m (\AA_\ast (\RR P(W)))  & \ar[l] \bL_n (\AA_\ast (\RR P(V\oplus
W))) }
\]
where $m=\dim(W)-1$ and $n=\dim(V\oplus W)-1$. It is also easy to
promote the resulting composite map
\[ \varphi\co \bS(\RR P(V\oplus W)) \ra \bL_m (\AA_\ast (\RR P(W))) \]
to a natural transformation between functors in the variable $V$
(note that the target functor is constant). More precisely we can
construct a continuous functor $G_W \co \sJ \ra \Spaces_\ast$ which
is polynomial of degree $0$ (that is, essentially constant) and a
natural transformation
\[ \varphi\co \bS(\RR P(V\oplus W)) \ra G_W(V) \]
which specializes to Ranicki's map $\bS(\RR P(W)) \ra \bL_m
(\AA_\ast (\RR P(W)))$ when $V=0$. Now we can state our Thom
isomorphism result.

\vspace{2mm} \noindent \textbf{Theorem B.} {\it Let $W$ in $\sJ$ be
oriented, of even dimension. There is a natural map
\[
\Omega^W F^a(V) \ra \hofiber[ F^a(V\oplus W) \stackrel{\varphi}{\ra}
G_W(V) ]
\]
which is a homotopy equivalence for $\dim(V)\ge 3$. }

\medskip\nin

The following is a simple reformulation of the main result of
\cite{Ma}.

\vspace{2mm} \noindent \textbf{Theorem C.} {\it Let $W \in \sJ$ be
such that $\dim (W) \geq 6$. Then the functor
\[
V \mapsto F^g (V \oplus W)
\]
on $\sJ$ is polynomial of degree $\leq 1$.} \vspace{1mm}

\begin{proof}[Proof of Theorem 1.] We assume a
natural homotopy fiber sequence \[ F^g(V)\to F^a(V)\to \ZZ\] when
$\dim(V)\ge 6$. We use that and Theorem C to deduce that for $W$ in
$\sJ$ with $\dim(W)\ge 6$, the functor $V \mapsto F^a(V \oplus W)$
is polynomial of degree $\leq 1$ without any low-dimensional
deviations. Now suppose in addition that $W$ is even-dimensional and
oriented. Let $F^a_W$ be  the functor taking $V$ in $\sJ$ to the
homotopy fiber of $\varphi\co F^a(V \oplus W) \rightarrow G_W(V)$ in
Theorem B. Then $F^a_W$ is polynomial of degree $\leq 1$ without any
low-dimensional deviations, because it is the homotopy fiber of a
natural transformation between a functor which is polynomial of
degree $\leq 1$ and another functor which is polynomial of degree
$\leq 0$. Therefore we have
\[
F^a_W (V) \simeq T_1 F^a_W (V)
\]
for all $V$ in $\sJ$. From Theorem B we obtain
\[
T_1 F^a_W (V) \simeq T_1 \Omega^W F^a(V)
\]
for all $V \in \sJ$, since $F^a_W$ and $\Omega^WF^a$ ``agree'' on
objects $V$ of sufficiently large dimension. Now suppose in addition
that $4$ divides $\dim (W)$. Then we get from Theorem A that
\[
 T_1 \Omega^W F^a(V) \simeq T_1 F^a(V)
\]
for all $V$ in $\sJ$. Composing these three natural homotopy
equivalences, we get $F^a_W(V) \simeq T_1F^a(V)$ for all $V$ in
$\sJ$. Specializing to $V=0$ and unraveling the definition of
$F^a_W(0)$ we obtain the statement of Theorem 1.
\end{proof}


\section{Overview of algebraic surgery}
\label{alg-sur}

The aim of this section is to recall the homotopy fibration sequence
of algebraic surgery due to Ranicki \cite{Ra}. It has the form
\begin{equation} \label{alg-sur-seq}
\bS (X,w) \rightarrow \bL_n (\AA_\ast (X,w) ) \rightarrow \bL_n
(\AA[\pi]^w),
\end{equation}
where $X$ is a connected simplicial complex (equipped with a
universal covering, with deck transformation group $\pi$) and $w
\colon \pi \rightarrow \ZZ_2$ is a homomorphism. More generally, $X$
can be a $\Delta$-complex (see \cite{Hat} and
definition~\ref{Deltacx} below). If $X$ is an $n$-dimensional
manifold and $w \colon \pi \rightarrow \ZZ_2$ is the orientation
character, then up to a small modification the homotopy fibration
sequence (\ref {alg-sur-seq}) can be identified with the geometric
homotopy fibration sequence of surgery
\begin{equation} \label{sur-seq}
\sStw (X) \rightarrow \sN (X) \rightarrow \sL_n (X).
\end{equation}

This section contains essentially no new results. It is a review of
definitions and tools we need from \cite{Ra} (see also \cite{We2}).
We focus mostly on the case where $w$ is trivial, but towards the
end of the section we indicate the modifications needed if $w$ is
not trivial (the non-orientable case). All spaces in
(\ref{alg-sur-seq}) are certain $L$-theory spaces associated with
various additive categories with chain duality and all the maps in
(\ref{alg-sur-seq}) are induced by functors between these
categories. We start with an additive category $\AA$ and introduce
the category $\BB(\AA)$ of chain complexes of $\AA$-objects, graded
over $\ZZ$ and bounded above and below. The notion of a chain
duality $T \colon \AA \rightarrow \BB(\AA)$ on $\AA$ will be
recalled below. We use this to define symmetric and quadratic
structures on objects of $\BB(\AA)$.

\begin{defn}
\label{defn-chduality} A {\it chain duality} on an additive category
$\AA$ is a contravariant additive functor $T \colon \AA \rightarrow
\BB (\AA)$ together with a natural transformation $e$ from $T^2
\colon \AA \rightarrow \BB (\AA)$ to $\id \colon \AA \rightarrow \BB
(\AA)$ such that for each $M$ in $\AA$
\begin{enumerate}
\item $e_{T(M)} \cdot T(e_M) = \id \colon T(M) \rightarrow T^3 (M) \rightarrow
T(M)$,
\item $e_M \colon T^2 (M) \rightarrow M$ is a chain homotopy equivalence.
\end{enumerate}
\end{defn}

A chain duality $T \colon \AA \rightarrow \BB(\AA)$ can be
canonically extended to the chain duality $T \colon \BB(\AA)
\rightarrow \BB(\AA)$ as follows. Let $C$ be a chain complex in
$\BB(\AA)$. Then we can define a double complex
\[
T(C)_{p,q} = T(C_{-p})_q.
\]
The dual chain complex $T(C) \in \BB(\AA)$ is the total complex of
this double complex.

A chain duality $T \colon \AA \rightarrow \BB(\AA)$ can be used to
define a tensor product of two objects $M$, $N$ in $\AA$ over $\AA$
as
\begin{equation} \label{tensor-product}
M \otimes_\AA N = \Hom_\AA (T(M),N).
\end{equation}
This is a chain complex of abelian groups. 

The main examples of additive categories with chain duality we will
consider are the following.

\begin{expl} \label{R-duality}
Let $R$ be a ring with involution $r\mapsto\bar r$ and let $\AA(R) $
be the category of finitely generated projective left $R$-modules.
On the category $\AA (R)$ we can define a chain duality $\AA\to \AA$
by $T(M) = \Hom_R (M,R)$. The involution can be used to make $T(M)$
into a f.g. projective left $R$-module. The dual $T(C)$ of a finite
chain complex $C$ in $\AA(R)$ is $\Hom_R(C,R)$.

The most important example for us is $R = \ZZ[\pi]$, the group ring
of a group $\pi$, with involution given by $\bar g=g^{-1}$ for
$g\in\pi$.

The category $\AA(\ZZ)$ with chain duality will sometimes be denoted
just by $\AA$, and the category $\AA(\ZZ[\pi])$ will sometimes be
denoted $\AA[\pi]$.
\end{expl}

In this paper we write $\Delta$ for the category with objects $\uli
n=\{0,1,\dots,n\}$, for $n=0,1,2,\dots$, and order-preserving
\emph{injective} maps as morphisms. A $\Delta$-set is a functor from
$\Delta\op$ to sets. A $\Delta$-set $Y$ has a geometric realization
$|Y|$. It is the quotient of $\coprod_n Y_n\times\Delta^n$ by the
relations $(u^*y,x)\sim (y,u_*x)$ for $y\in Y_n$~, $x\in \Delta^m$
and $u\co\uli m\to\uli n$ a morphism in $\Delta$. \newline Out of a
$\Delta$-set $Y$, we can make a category $\cat(Y)$ with object set
$\coprod_n Y_n$~, where a morphism from $\sigma\in Y_m$ to $\tau\in
Y_n$ is a morphism $u\co \uli m\to \uli n$ in $\Delta$ with
$f^*\tau=\sigma$. We write $u\colon \sigma\to \tau$
for short. 

\begin{defn}
\label{Deltacx}
 A $\Delta$-complex is a space $X$ together with a $\Delta$-set $sX$
and a homeomorphism $|sX|\to X$. It is considered \emph{finite} if
$sX$ is finite. When we write \emph{simplex in X}, for a
$\Delta$-complex $X$, we mean a simplex in $sX$.
\end{defn}

\begin{expl} \label{X-based-duality}
Let $\AA$ be an additive category and let $X$ be a finite
$\Delta$-complex. Then there are defined two additive categories
$\AA_\ast(X)$ and $\AA^\ast(X)$ of $X$-based objects in $\AA$. An
object $M$ of $\AA$ is $X$-{\it based} if it comes as
\[
M = \sum_{n\ge 0}\sum_{\,\,\sigma \in sX_n} M(\sigma).
\]
A morphism $f\co M\rightarrow N$ in $\AA^\ast (X)$, resp.\ $\AA_\ast
(X)$, is a matrix $f = (f_u)$ of morphisms $f_u\co M(\sigma)
\rightarrow N(\tau)$ in $\AA$, resp. $f_u\co M(\tau)\rightarrow
N(\sigma)$ in $\AA$, with entries corresponding to morphisms $u\co
\sigma\to \tau$ between simplices of $X$. Composition of morphisms
is given by matrix multiplication.
\newline
Such a morphism $f$ can be thought of as an upper triangular, resp.\
lower triangular matrix. For example, $f$ is an isomorphism if and
only if all diagonal entries, the $f_u$ in which $u$ is an identity,
are invertible in $\AA$. \newline Given $N$ in $\AA$ and $\sigma$ in
$X$, let $N_\sigma$ in $\AA_*(X)$, resp. $\AA^*(X)$, be defined by
$N_\sigma(\sigma)=N$ and $N_\sigma(\tau)=0$ for $\tau\ne\sigma$.
Clearly $N\to N_\sigma$ is a functor from $\AA$ to $\AA_*(X)$, resp.
$\AA^*(X)$. This functor has a right adjoint $M\to M[\sigma]$ from
$\AA_*(X)$, resp. $\AA^*(X)$, to $\AA$. We have
$M[\sigma]=\sum_{\sigma\to\tau}M(\tau)$~, resp.
$M[\sigma]=\sum_{\tau\to\sigma}M(\tau)$ where the direct sum is
taken over all morphisms $\sigma\to\tau$, resp. $\tau\to\sigma$,
with fixed $\sigma$ and arbitrary $\tau$. For a morphism $f\co
M_0\to M_1$ in $\AA_*(X)$, resp. $\AA^*(X)$, the induced morphism
$M_0[\sigma] \to M_1[\sigma]$ is a sum of terms $f_u\colon
M_0(\tau)\to M_1(\rho)$, one such for every diagram
\[ \sigma \to \tau \stackrel{u}{\lra} \rho~, \textup{ resp. }\,\,
\rho \stackrel{u}{\lra} \tau\to\sigma \] of simplices in $X$.
\newline Rather than defining chain dualities in $\AA_\ast(X)$ and
$\AA^\ast(X)$ directly, we focus on the tensor products,
$\otimes_{\AA_\ast(X)}$ and $\otimes_{\AA^\ast(X)}$, which are
easier to motivate. Suppose therefore that $M$ and $N$ are $X$-based
objects of $\AA$. Then
\[
\begin{array}{ccc}
(M \otimes_{\AA_\ast(X)} N)_r & = & \sum_{\sigma \in sX_r}
\sum_{\lambda\leftarrow \sigma\to\mu} (M(\lambda) \otimes_\AA
N(\mu))_{r-|\sigma|} \\
(M \otimes_{\AA^\ast(X)} N)_r & = &\prod_{\sigma \in sX_r}
\sum_{\lambda\to \sigma\leftarrow \mu} (M(\lambda) \otimes_\AA
N(\mu))_{r+|\sigma|}.
\end{array}
\]
In the case where the $\Delta$-complex $X$ is a simplicial complex
(with ordered vertex set, say), these graded abelian groups can be
regarded as chain subcomplexes of
$C_*X\otimes_{\ZZ}(M\otimes_{\AA}N)$ and
$\Hom_{\ZZ}(C_*X,M\otimes_{\AA}N)$ respectively, where $C_*X$ is the
cellular chain complex of $X$. In the general case, we can still say
that
\[
\begin{array}{ccc}
\sigma \mapsto \sum_{\lambda\leftarrow \sigma\to\mu} M(\lambda)\otimes_\AA N(\mu) \\
\sigma \mapsto \sum_{\lambda\to \sigma\leftarrow \mu} M(\lambda)
\otimes_\AA N(\mu)
\end{array}
\]
is a contravariant (resp. covariant) functor, with chain complex
values, on the category of simplices of $X$. This determines in the
usual way a double chain complex of abelian groups. The
corresponding total chain complex is $M \otimes_{\AA_\ast (X)} N$,
resp. $M \otimes_{\AA^\ast (X)} N$. The adjunction
(\ref{tensor-product}) then determines the chain duality functors
$\AA_\ast(X)\to \BB(\AA_\ast(X))$ and $\AA^\ast(X)\to
\BB(\AA^\ast(X))$ as follows. Let $M$ be an object in $\AA_\ast(X)$,
resp.\ $\AA^\ast(X)$. Then
\[
T(M)_r (\sigma) = \begin{cases} T(M[\sigma])_{r+|\sigma|},\\
T(M[\sigma])_{r-|\sigma|}, \end{cases}
\]
with differential
\[ d_{T(M)}(u\colon\sigma\to\tau):\quad
\begin{cases} T(M[\tau]) \to T(M[\sigma]) \\
T(M[\sigma]) \to T(M[\tau])
\end{cases}
\]
equal to $d_{T(M[\sigma])}$ if $\sigma=\tau$, equal to
\[
\begin{cases}
(-1)^iT\big(u^*\colon M[\tau] \to M[\sigma]\big) \\
(-1)^iT\big(u_*\colon M[\sigma] \to M[\tau]\big)
\end{cases}
\]
if $|\tau|=|\sigma|+1$ and $u$ omits the $i$-th vertex, and equal to
0 for all other $\sigma\to\tau$.
\end{expl}

\begin{rem} \label{cofibrant}
An object $M$ in $\AA_*(X)$ determines a contravariant functor
$M_\natural$ from the category of simplices of $X$ to $\AA$ by
$M_\natural(\sigma) = M[\sigma]$. We call such a contravariant
functor from the category of simplices of $X$ to $\AA$, or any
isomorphic one, \emph{cofibrant}. Similarly an object $M$ in
$\AA^*(X)$ determines a covariant functor $M^{\natural}$ from the
category of simplices of $X$ to $\AA$ by $M^\natural(\sigma) =
M[\sigma]$. Again we call such a covariant functor, or any
isomorphic one, cofibrant. A morphism $f\colon M\to N$ in
$\AA_\ast(X)$ induces a natural transformation $f_\natural\co
M_\natural\to N_\natural$, and vice versa. A morphism $f\co M\to N$
in $\AA^\ast(X)$ induces a natural transformation $f^\natural\co
M^\natural\to N^\natural$, and vice versa. In this way $\AA_\ast(X)$
and $\AA^\ast(X)$ are equivalent to, and could be re--defined as,
certain categories of functors on the category of simplices of $X$.
There are situations when we have to resort to these alternative
definitions.
\end{rem}

\begin{expl} \label{diss-duality}
For more motivation of the duality on the categories of $X$-based
objects here is an example. To start with let $X$ be a finite
simplicial complex with ordered vertex set and let $X'$ be its
barycentric subdivision. The simplices $\sigma$ of $X$ correspond to
the vertices $\hat\sigma$ of $X'$. For a simplex $\sigma$ of $X$ its
dual cell $D(\sigma,X)$ is the subcomplex of $X'$ spanned by the
simplices with vertex set of the form
\[ \{\hat\tau_0,\hat\tau_1,\dots,\hat\tau_p\} \]
where $\tau_0$ contains $\sigma$ and
$\tau_0\subset\tau_1\subset\cdots\subset \tau_p$. Its ``boundary''
is spanned by all simplices in $D(\sigma,X)$ which do not have
$\hat\sigma$ as a vertex. The dual cell $D(\sigma,X)$ is
contractible. Apart from that it does not always have the properties
that we would expect from a cell (such as being homeomorphic to a
euclidian space), but it has the dual properties. In particular,
$D(\sigma,X)\smin\partial D(\sigma,X)$ has a trivial normal bundle
in $X$ with fibers homeomorphic to $\RR^{|\sigma|}$.

Let $M$ be a closed $n$-dimensional smooth or $PL$ manifold. Any map
$f \colon M \rightarrow X$ is homotopic to a map transverse to the
dual cells of $X$. If $f$ is transverse to the dual cells, then
\[
(M[\sigma],\partial M[\sigma]) = f^{-1} (D(\sigma,X),\partial
D(\sigma,X))
\]
is an $(n-|\sigma|)$-dimensional manifold. The collection $\{
M[\sigma] \; | \; \sigma \textup{ in } X \}$, or more precisely, the
contravariant functor $\sigma\mapsto M[\sigma]$, is called an
$X$-dissection of $M$. In this situation there exists a structure of
a $CW$-space on $M$ such that each $M[\sigma]$ is a $CW$-subspace.
The cellular chain complex $C_*M$ can then be understood as a chain
complex in $\BB (\AA_\ast (X))$ via the decomposition
\[
C_*M = \sum_{\sigma} C_*(M[\sigma],\partial M[\sigma]).
\]
We may expect this to be self-dual, with a shift of $n$, since $M$
is a closed manifold. The dual of $C_*M$ in $\AA_\ast (X)$ is by
definition
\[
T(C_*M) = \sum_{\sigma} C^{-|\sigma|-\ast}(M[\sigma]).
\]
The ordinary Poincar\'e duality homotopy equivalences
\[ C_*(M[\sigma],\partial M[\sigma])\simeq C^{n-|\sigma|-\ast}(M[\sigma]) \]
suggest that $\Sigma^nT(C_*M)$ is indeed homotopy equivalent in
$\AA_\ast (X)$ to $C_*M$. This will be confirmed later. \newline Now
we need to generalize these observations from the setting of
simplicial complexes to that of $\Delta$-complexes. For a
$\Delta$-complex $X$ and a simplex $\sigma$ in $X$, we have the
category of simplices of $X$ under $\sigma$. Its objects are
morphisms $u\co\sigma\to \tau$ where $\sigma$ is fixed and $\tau$ in
$X$ is variable. Its nerve is a $\Delta$-set and the corresponding
$\Delta$-complex is, by definition, the dual cell $D(\sigma,X)$. The
boundary $\partial D(\sigma,X)$ corresponds to the nerve of the full
subcategory with objects $u\colon\sigma\to \tau$ where $u$ is not an
identity. \newline The dual cell $D(\sigma,X)$ is contractible,
because the category of simplices of $X$ under $\sigma$ has an
initial object. There is a canonical map
\[  c_\sigma\colon D(\sigma,X) \lra X \]
defined as follows. A $k$-simplex of $D(\sigma,X)$ corresponds to a
diagram
\[  \sigma\to \tau_0\to\tau_1\to \cdots\to\tau_k \]
of simplices in $X$. The vertices of that $k$-simplex are the
resulting $\sigma\to\tau_i$ for $i=0,1,\dots,k$. The restriction of
$c_\sigma$ to the $k$-simplex is the ``linear'' map taking the
vertex $\sigma\to\tau_i$ to the barycenter of $\tau_i$ in $X$.
\newline The map $c_\sigma$ need not be injective. However, it is
locally injective, it embeds $D(\sigma,X)\smin \partial
D(\sigma,X)$, and the image of that restricted embedding has a
trivialized normal bundle in $X$, with fibers homeomorphic to
$\RR^{|\sigma|}$. This results in a stratification of $X$ where the
strata have the form
\[ X(\sigma)= c_\sigma(D(\sigma,X)\smin \partial D(\sigma,X)) \]
and each stratum has a trivialized normal bundle with fiber
homeomorphic to a euclidian space. The closure $X[\sigma]$ of
$X(\sigma)$ in $X$ is the union of all $X(\tau)$ for which there
exists a morphism $\sigma\to \tau$. (Zeeman's dunce hat, the
two-dimensional $\Delta$-complex with a single 0-simplex, a single
1-simplex and a single 2-simplex, is an instructive example.)
\newline Let $M$ be an $n$-dimensional smooth or $PL$ manifold. Any
map $f \colon M \rightarrow X$ is homotopic to a map transverse to
the stratification of $X$ by subsets $X(\sigma)$. If $f$ is
transverse to the stratification, then the pullback of $f$ and
$c_\sigma$ is a manifold $M[\sigma]$ with boundary $\partial
M[\sigma]$. Any morphism $\sigma\to \tau$ in $X$ determines a map
$M[\tau]\to M[\sigma]$ which, if $|\tau|>|\sigma|$, factors through
$\partial M[\sigma]$. That map is locally an embedding and it embeds
$M[\tau]\smin\partial M[\tau]$. The functor $\sigma\mapsto
M[\sigma]$ together with the identification $\colim_{\sigma}
M[\sigma]\cong M$ is called an $X$-dissection of $M$. It is always
possible to equip the functor $\sigma\mapsto M[\sigma]$ with a
CW-structure. (A CW-structure on a contravariant functor $F$ from a
small category to spaces is a filtration of $F$ by subfunctors $F_i$
for $i=-1,0,1,2,3,\dots$, where $F_{-1}=\emptyset$ and $F_i$ is
obtained from $F_{i-1}$ by ``attaching'' functors of the form
$a\mapsto D^i\times \coprod_{\lambda}\hom(a,b_{\lambda})$ using
natural attaching maps $S^i\times
\coprod_{\lambda}\hom(a,b_{\lambda})\to F_{i-1}(a)$. If $F$ comes
with a CW-structure, we also say that $F$ is a CW-functor. See
\cite{Dro} for more details.) The cellular chain complex $C_*M$ can
then be understood as a chain complex in $\BB (\AA_\ast(X))$ via the
decomposition
\[
C_*M = \sum_{\sigma} C_*(M[\sigma],\partial M[\sigma]).
\]
where $C_*(M[\sigma],\partial M[\sigma])$ is the cellular chain
complex of $M[\sigma]/\partial M[\sigma]$. The dual of $C_*M$ in
$\AA_\ast(X)$ is
\[
T(C_*M) = \sum_{\sigma} C^{-|\sigma|-\ast}(M[\sigma]).
\]
If $M$ is (only) a topological manifold, and $f\co M\to X$ is
transverse to the stratification of $X$ by subsets $X(\sigma)$, then
we can still construct a contravariant functor $\sigma\mapsto
F[\sigma]$ with CW-structure from the category of simplices of $X$
to the category of spaces, and a natural homotopy equivalence
$F[\sigma]\to M[\sigma]$, for $\sigma$ in $X$. Then $\partial
F[\sigma]$ is well defined: it is the CW-subspace of $F[\sigma]$
containing all the cells which come from some $F[\tau]$ via some
$u\colon\sigma\to\tau$. The object
\[
C_*F = \sum_{\sigma} C_*(F[\sigma],\partial F[\sigma])
\]
in $\BB(\AA_\ast(X))$ is a good substitute for a possibly
nonexistent $C_*M$.
\end{expl}

\medskip
A chain duality $T \colon \AA \rightarrow \BB (\AA)$ can be used to
define symmetric and quadratic chain complexes in $\BB(\AA)$ as
follows. Firstly, notice that given two objects $M$ and $N$ of
$\AA$, their tensor product $M \otimes_\AA N$ possesses a symmetry
isomorphism
\[
T_{M,N} \colon M \otimes_{\AA} N \rightarrow N \otimes_{\AA} M
\]
given by taking
\[
f \in (M \otimes_{\AA} N)_n = \Hom_{\AA} (T(M)_{-n},N)
\]
to
\[
T_{M,N} (f) \in (N \otimes_{\AA} M)_n = \Hom_{\AA} (T(N)_{-n},M)
\]
where
\[
T_{M,N} (f) = e_M \cdot T(f) \colon T(N)_n \rightarrow
T(T(M)_{-n})_{-n} \subseteq T^2 (M)_0 \rightarrow M.
\]
This tensor products extends to a tensor product $C \otimes_{\AA} D$
of chain complexes $C,D$ in $\BB (\AA)$ and there is also a symmetry
isomorphism
\[
T_{C,D} \colon C \otimes_\AA D \rightarrow D \otimes_\AA C.
\]
If $C=D$ this makes $C \otimes_{\AA} C$ into a finite chain complex
of $\ZZ[\ZZ_2]$-modules. Now let $W$ be the standard
$\ZZ[\ZZ_2]$-resolution of $\ZZ$, i.e. it is a chain complex of
$\ZZ[\ZZ_2]$-modules
\[
W = \cdots \ZZ[\ZZ_2] \xrightarrow{1+T} \ZZ[\ZZ_2] \xrightarrow{1-T}
\ZZ[\ZZ_2] \rightarrow 0
\]
concentrated in non-negative degrees. Then there are the following
two chain complexes of abelian groups
\[
\Hom_{\ZZ[\ZZ_2]} (W, C \otimes_{\AA} C),
\]
\[
W \otimes_{\ZZ[\ZZ_2]} (C \otimes_{\AA} C).
\]

\begin{defn}
\label{defn-SAPC} An $n$-dimensional {\it symmetric algebraic
complex} in $\BB(\AA)$ is a pair $(C,\phi)$ with $C$ a chain complex
in $\BB(\AA)$ and $\phi$ an $n$-cycle in $\Hom_{\ZZ[\ZZ_2]} (W, C
\otimes_{\AA} C)$. An $n$-dimensional {\it quadratic algebraic
complex} in $\BB(\AA)$ is a pair $(C,\psi)$ with $C$ a chain complex
in $\BB(\AA)$ and $\psi$ an $n$-cycle in $W \otimes_{\ZZ[\ZZ_2]} (C
\otimes_{\AA} C)$.
\end{defn}

We note that in the above definition it is not required that the
chain complex $C$ is concentrated in dimensions from $0$ to $n$, it
is only required that it is bounded below and above. The dimension
$n$ is associated with the symmetric structure $\phi$ or with the
quadratic structure $\psi$.

An $n$-dimensional symmetric structure $\phi$ on a chain complex $C$
can be described as a collection of chains in $\Hom_{\AA}(T(C),C)$,
\[
\phi = \{\phi_s \colon T(C)_{-\ast} \rightarrow C_{n-\ast+s} \; |
\;s \geq 0 \}
\]
satisfying certain relations.

An $n$-dimensional quadratic structure $\psi$ on a chain complex $C$
can be described as a collection of chains in $\Hom_{\AA}(T(C),C)$,
\[
\psi = \{\psi_s \colon T(C)_{-\ast} \rightarrow C_{n-\ast-s} \; |
\;s \geq 0 \}
\]
satisfying certain relations.

An $n$-dimensional quadratic structure $\psi$ on $C$ determines an
$n$-dimensional symmetric structure $\phi$ on $C$ by
$\phi_0=(1+T)\psi_0$ and $\phi_s=0$ for $s>0$. We describe this
relationship by writing $\phi=(1+T)\psi$.

\begin{defn} For $C$ in $\BB(\AA)$, an $n$-cycle in $C \otimes_{\AA}C
\cong \Hom_{\AA}(TC,C)$ is \emph{nondegenerate} if the corresponding
chain map $TC\to C$ of degree $n$ is a chain homotopy equivalence.
An $n$-dimensional {\it symmetric algebraic Poincar\'e complex}
(SAPC) in $\BB(\AA)$ is a symmetric algebraic complex $(C,\phi)$
such that $\phi_0$ is nondegenerate. An $n$-dimensional {\it
quadratic algebraic Poincar\'e complex} (QAPC) in $\BB(\AA)$ is a
quadratic algebraic complex $(C,\psi)$ such that $(1+T)\psi_0$ is
nondegenerate.
\end{defn}

\begin{expl} \label{sym-con}
Let $X$ be a connected finite $CW$-complex and $\widetilde X\to X$ a
universal covering with deck transformation group $\pi$. The
diagonal map
\[ \nabla\co X\to \widetilde X\times_{\pi}\widetilde X \]
is a $\ZZ_2$-map for the trivial action of $\ZZ_2$ on the source and
the permutation action on the target. It is not cellular in general.
However it is easy to construct a \emph{cellular} $\ZZ_2$-map
\[
\nabla^{\sharp}\co E\ZZ_2 \times X \ra \widetilde
X\times_{\pi}\widetilde X \] which is $\ZZ_2$-homotopic to the
composition of $\nabla$ just above with the projection $ E\ZZ_2
\times X\ra X$. Here $E\ZZ_2$ can be taken as the universal
(=double) cover of $B\ZZ_2=\RR P^{\infty}$, with the standard
$CW$-structure on $\RR P^{\infty}$. Hence the map of cellular chain
complexes induced by $\nabla^{\sharp}$ takes the form
\[ W\otimes C_*X  \lra C_*\widetilde X\otimes_{\ZZ[\pi]}C_*\widetilde X\]
with adjoint
\[ C_*X\lra \Hom_{\ZZ[\ZZ_2]}(W,
C_*\widetilde X\otimes_{\ZZ[\pi]}C_*\widetilde X).
\]
Regard now $C_* \widetilde X$ as an object in $\BB(\AA[\pi])$, with
$\AA[\pi]$ as in example~\ref{R-duality}. Then by all the above, any
$n$-cycle $\mu$ in $C_*X$ determines an $n$-dimensional symmetric
structure $\phi (X)$ on $C_* \widetilde X$. If $X$ is an orientable
Poincar\ee duality space and $\mu$ represents a fundamental class
$[X]$, then $\phi_0$ is nondegenerate and so
\[  (C_* \widetilde X,\phi(X))\]
is an $n$-dimensional SAPC.
\end{expl}

\begin{expl} \label{quad-con}
Let $(f,b) \colon M \rightarrow X$ be a degree one normal map of
$n$-dimensional closed manifolds or Poincar\ee duality spaces, where
$X$ is connected and equipped with a universal covering. Denote by
$K(f)$ the algebraic mapping cone of the Umkehr map of chain
complexes
\[
\CD f^! \colon C_\ast\widetilde{X} \simeq C^{n - \ast} \widetilde{X}
@>\quad{f^{n -\ast}}\quad>> C^{n - \ast} \widetilde{M} \simeq C_\ast
\widetilde{M}.
\endCD
\]
As explained just above, $C_*\widetilde X$ comes with a structure of
$n$-dimensional SAPC over $\ZZ[\pi]$. This projects to a structure
of $n$-dimensional SAPC on $K(f)$. Ranicki in \cite{RaLMS2},
\cite{Ra} refines the latter to an $n$-dimensional QAPC on
$(K(f),\psi(f))$.
\end{expl}

In the next definition, the standard simplex $\Delta^n$ is regarded
as a simplicial complex in the usual way. Each face inclusion
$\Delta^{n-1}\to \Delta^n$ induced by the monotone injection
$\{0,1,\dots,n-1\}\to \{0,1,\dots,n\}$ induces an additive functor
$d_i\co \AA^\ast(\Delta^n)\to \AA^\ast(\Delta^{n-1})$ which commutes
with the chain dualities. We identify $\AA^\ast(\Delta^0)$ with
$\AA$.

\begin{defn}
Two $n$-dimensional SAPC (QAPC) in $\BB(\AA)$, say $(C,\phi)$ and
$(C',\phi')$, are called {\it cobordant} if there exists an
$n$-dimensional SAPC (QAPC), say $(D,\psi)$ in
$\BB(\AA^\ast(\Delta^1))$, such that $d_0(D,\psi) \cong (C,\phi)$
and $d_1 (D,\psi)\cong (C',\phi')$.
\end{defn}

With the alternative definition of $\AA^\ast(\Delta^1)$ outlined in
remark~\ref{cofibrant}, the cobordism relation is an equivalence
relation on $n$-dimensional SAPC (QAPC). The direct sum makes the
cobordism classes of SAPC (QAPC) into an abelian group, where the
inverse of $[(C,\phi)]$ is given by $[(C,-\phi)]$.

\begin{defn}
The group of cobordism classes of $n$-dimensional SAPC in $\BB(\AA)$
is denoted by $L^n (\AA)$. The group of cobordism classes of
$n$-dimensional QAPC in $\BB(\AA)$ is denoted by $L_n (\AA)$.
\end{defn}

\begin{expl}
For the category $\AA[\pi]$ with chain duality as in Example
\ref{R-duality}, the $L$-groups $L^n (\AA[\pi])$ are the usual
symmetric $L$-groups $L^n(\pi)$ of Mishchenko and the $L$-groups
$L_n (\AA[\pi])$ are the quadratic $L$-groups $L_n(\pi)$ of Wall
(see \cite{Ra}).
\end{expl}

The $L$-groups are in fact the homotopy groups of certain spaces.
These are defined as $\Delta$-sets (alias ``simplicial sets without
degeneracy operators'') as follows.

\begin{defn}
Let $\bL^n (\AA)$, resp. $\bL_n (\AA)$, denote the $\Delta$-set
whose $k$-simplices are $n$-dimensional SAPC, resp. QAPC in the
category $\AA^\ast (\Delta^k)$. The face maps are induced by the
functors $d_i \colon \AA^\ast (\Delta^k) \rightarrow \AA^\ast
(\Delta^{k-1})$. We use the alternative definitions of $\AA^\ast
(\Delta^k)$ given in remark~\ref{cofibrant}.
\end{defn}

\begin{expl}
With $\AA[\pi]$ as in Example \ref{R-duality}, the $L$-theory space
$\bL_n (\AA[\pi])$ is the $L$-theory space $\bL_n (\pi)$ of Quinn,
with homotopy groups $\pi_k\bL_n (\pi)=L_{k+n} (\pi)$.
\end{expl}

\begin{rem} \label{quad-sym}
The assignment $\psi \mapsto (1+T) \cdot \psi$ for $(C, \psi)$ an
$n$-dimensional QAPC in $\BB (\AA)$ defines a symmetrization map
$(1+T) \colon \bL_n (\AA) \rightarrow \bL^n (\AA)$.
\end{rem}

Let $X$ be a finite Delta-complex. We now have the definitions of
the spaces $\bL_n(\AA_\ast(X))$ and $\bL_n(\AA[\pi])$ from the
homotopy fibration sequence (\ref{alg-sur-seq}), ignoring
orientation matters which will be discussed later. Assuming
$\AA=\AA(\ZZ)$ for simplicity, we proceed to describe the map
$\alpha$ from $\bL_n(\AA_\ast (X))$ to $\bL_n(\AA[\pi])$ which is
called {\it assembly}.

Suppose that $X$ comes with a principal $\pi$--bundle $p\co
X^\natural \to X$. (In most applications this will be a universal
covering for $X$, and $X$ will be connected, but we do not have to
assume that now.) The map $\alpha$ is induced by an additive
functor, also denoted  $\alpha$. Define
\[
\alpha \colon \AA_\ast (X) \rightarrow \AA [\pi]
\]
on objects by
\[\alpha(M)= \sum_{\sigma\in sX^\natural}
M(p(\sigma)) \] with $\pi$ acting on the right-hand side by
permutating summands in the obvious way. A morphism $f\co M\to N$ in
$\AA_\ast(X)$ induces $\alpha(f) \colon \alpha(M) \rightarrow
\alpha(N)$, which we define in matrix notation by
\[\alpha(f)_{(\sigma,\tau)}= \sum_{u\colon\sigma\to\tau}f_{p(u)} \]
where $u\colon \sigma\to\tau$ in $X^\natural$ and $p(u)\colon
p(\sigma)\to p(\tau)$ is the induced morphism in $X$. \newline In
order to see that the assembly functor $\alpha$ induces a map of the
$L$-spaces one has to see that it ``commutes'' with the chain
dualities as in Examples \ref{R-duality}, \ref{X-based-duality}, or
with the corresponding tensor products. We choose the tensor product
option. The coefficient system $\sigma\mapsto
\sum_{\lambda\leftarrow\sigma\to\mu}M(\lambda)\otimes_{\AA}N(\mu)$
on $X$ comes with an evident natural transformation to the constant
coefficient system $\sigma\mapsto
\alpha(M)\otimes_{\AA[\pi]}\alpha(N)$, due to the fact that a
diagram such as $\lambda\leftarrow\sigma\to\mu$ in $X$ determines a
preferred path class in $X$ connecting the barycenters of $\lambda$
and $\mu$. Passing to the cellular chain complexes associated with
these coefficient systems gives
\[ M\otimes_{\AA_\ast(X)}N \lra
C_*X\otimes\big(\alpha(M)\otimes_{\AA[\pi]}\alpha(N)\big). \] We
compose with the augmentation $C_*(X)\to \ZZ$ to get
\[ M\otimes_{\AA_\ast(X)}N \lra \alpha(M)\otimes_{\AA[\pi]}\alpha(N) \]
and more generally
\[ C\otimes_{\AA_\ast(X)}D \lra \alpha(C)\otimes_{\AA[\pi]}\alpha(D) \]
for objects $C$ and $D$ in $\BB(\AA)$. We use this to transport
symmetric and quadratic structures. Nondegeneracy is preserved, so
that QAPC are mapped to QAPC. It follows that we have a well defined
map of $L$-spaces
\begin{equation} \label{asmb}
\alpha \colon \bL_n (\AA_\ast(X)) \rightarrow \bL_n(\AA[\pi]),
\end{equation}
which is also called assembly. It is an algebraic version, due to
Ranicki, of the assembly map of Quinn \cite{Qu}. Apart from being
algebraic, it also incorporates Poincar\'e duality to switch from a
cohomological setup to a homological one.

\begin{rem}
In the above construction we indicated how an $n$-dimensional QAPC
in $\BB(\AA_\ast (X))$, say $(C,\psi)$, determines an assembled
$n$-dimensional QAPC in $\BB(\AA[\pi])$, denoted $\alpha (C,\psi)$.
For the sake of readability we will sometimes omit the prefix
$\alpha$ in the sequel, provided it is clear enough in which
category we are working.
\end{rem}

\begin{rem} \label{hlgy}
The spaces $\bL_n(\AA)$ can be arranged into an $\Omega$-spectrum
$\bL_{\bullet}(\AA)$, with homotopy groups $\pi_k\bL_{\bullet}(\AA)
\cong L_k(\AA)$ for $k\in \ZZ$. Beware that $\bL_n(\AA)$ is the
$(-n)$-th space in the $\Omega$-spectrum $\bL_{\bullet}(\AA)$. With
our conventions, $L_k(\AA)$ is isomorphic to $L_{k+4}(\AA)$ for all
$k\in \ZZ$, and indeed $\bL_{\bullet}(\AA)\simeq
\Omega^4\bL_{\bullet}(\AA)$. We also have Ranicki's law
\[
\pi_k (\bL_n(\AA_\ast(X))) \cong H_{n+k} (X;\bL_{\bullet}(\AA))
\]
for $k,n\in\ZZ$, and to be more precise $\bL_{\bullet}(\AA_\ast(X))
\simeq X_+\wedge \bL_{\bullet}(\AA)$. See \cite{Ra} for details.
When $\AA$ is the category of f.g. free $\ZZ$-modules with the
standard chain duality, then we write $\bL_{\bullet}$ for
$\bL_{\bullet}(\AA)$.
\end{rem}

\begin{expl} \label{frag-quad-con}
Let $(f,b) \colon M \rightarrow N$ be a degree one normal map of
$n$-dimensional closed smooth or PL manifolds and $g\colon N
\rightarrow X$ be a map to a simplicial complex $X$ such that both
$gf$ and $g$ are transverse to the dual cells of $X$. By Example
\ref{diss-duality} we have $X$-dissections $M\cong \colim M[\sigma]$
and $N\cong \colim N[\sigma]$, so that $C_*M$ and $C_*N$ can be
regarded as objects in $\BB(\AA_\ast(X))$, for suitable
$CW$--structures on $M$ and $N$. By analogy with
Example~\ref{sym-con}, there are preferred structures of
$n$--dimensional SAPC on $C_*M$ and $C_*N$, as objects of
$\BB(\AA_\ast(X))$. By analogy with Example~\ref{quad-con}, there is
an algebraic Umkehr map
\[ f^!\co C_*N \lra C_*M \]
in $\BB(\AA_\ast(X))$ with mapping cone $K(f)$, say. The resulting
structure of $n$--dimensional SAPC on $K(f)$, as an object of
$\BB(\AA_\ast(X))$, has a preferred refinement to a QAPC structure
$\psi$. We remark that $K(f)(\sigma)$ for a simplex $\sigma$ in $X$
can be identified with the mapping cone of an algebraic Umkehr map
\[ C_*(N[\sigma],\partial N[\sigma]) \lra C_*(M[\sigma],\partial M[\sigma])\]
which is the diagonal entry $f^!(\sigma,\sigma)$ of $f^!$. See again
\cite{Ra} for details. Under assembly, these constructions match and
recover those in Examples~\ref{sym-con} and~\ref{quad-con}.
\end{expl}

Suppose that $X$ is a connected $\Delta$-complex equipped with a
universal covering, with deck transformation group $\pi$. One way to
define the space $\bS(X,n)$ is to say that it is the homotopy fiber
of the assembly map
\[
\alpha\co \bL_n(\AA_\ast(X))\to \bL_n(\AA[\pi]).
\]
However, the theory of algebraic bordism categories of \cite{Ra} can
be used to provide a more direct description of the space $\bS(X,n)$
as the $L$-theory space associated to certain additive category with
chain duality (with certain restrictions on the objects). The
description is as follows.

\begin{defn}
\label{defn-struc} A $k$-simplex of the space $\bS(X,n)$ is an
$n$-dimensional QAPC in the category $\BB ((\AA_\ast (X))^\ast
(\Delta^k))$ which assembles to a contractible QAPC in
$\BB((\AA[\pi])^\ast(\Delta^k))$.
\end{defn}

Then there is an obvious inclusion map $\bS(X,n) \rightarrow \bL_n
(\AA_\ast(X))$. By the result of \cite[Proposition 3.9]{Ra} the
sequence (\ref{alg-sur-seq}) consisting of this map and the assembly
map $\alpha$ is a homotopy fibration sequence.

\medskip
Theorem 4.5 of \cite{We2} provides the following alternative to
definition~\ref{defn-struc}:

\begin{defn}
\label{defn-otherstruc} A $k$-simplex of $\bS(X,n)$ is an
$n$-dimensional SAPC in the category $\BB ((\AA_\ast (X))^\ast
(\Delta^k))$ which assembles to a contractible SAPC in
$\BB((\AA[\pi])^\ast(\Delta^k))$.
\end{defn}

\

\nin\textbf{Twisted versions.} Now we recall modifications in the
above machinery needed to treat the general case of nonorientable or
just nonoriented manifolds. It will be necessary to modify the
definition of the tensor product of $X$-based objects and the tensor
product of $\ZZ[\pi]$-modules and thus also the assembly map.

\begin{defn} \label{twist} A \emph{twist} on a group $\pi$ is a
$\pi$-module $\Gamma$ whose underlying abelian group is infinite
cyclic. A homomorphism of twisted groups, say from $(\pi,\Gamma)$ to
$(\pi',\Gamma')$, is a homomorphism $f\colon \pi\to \pi'$ together
with an isomorphism $\Gamma\to f^*\Gamma'$ of $\pi$-modules.
\end{defn}

\begin{expl} \label{loaded-manifold} For a
connected $n$--manifold $X$ and a universal covering of $X$ with
deck transformation group $\pi$, there is a canonical way to define
a twist on $\pi$. Let $\Gamma$ be the $n$--th integer homology with
locally finite coefficients of the universal covering. The action of
$\pi$ on $\Gamma$ is obvious.
\end{expl}

We now fix a finite $\Delta$-complex $X$ with a principal
$\pi$-bundle $p\colon X^\natural\to X$ and a twist $\Gamma$ on
$\pi$. The twist determines a homomorphism $w\colon\pi\to\{\pm1\}$
such that $gz=w(g)\cdot z$ for all $g\in\pi$ and $z\in\Gamma$. On
the group ring $\ZZ[\pi]$ we have the $w$-twisted involution given
by
\[
g \mapsto w(g) \cdot g^{-1}.
\]
The group ring with this involution will be denoted by $\ZZ[\pi]^w$
and the category of f.g. free left $\ZZ[\pi]^w$-modules will be
denoted by $\AA[\pi]^w$. We already have a chain duality $T$ on
$\AA[\pi]^w$ from Example~\ref{R-duality}; unfortunately this is no
longer considered quite right and to correct it we compose with the
functor $\Gamma\otimes_{\ZZ}$. Note that this causes a small change
in $\otimes_{\AA[\pi]^w}$ as well. Namely, for objects $M$ and $N$
in $\AA[\pi]^w$ there is the following isomorphism of abelian
groups:
\begin{equation} \label{iso}
M \otimes_{\AA[\pi]^w} N = \Gamma\otimes_{\ZZ[\pi]} (M
\otimes_{\ZZ}N).
\end{equation}

Similar remarks apply to $\AA_\ast(X)$, which we now rename
$\AA_\ast(X,w)$ to indicate a modified chain duality. We already
have a chain duality $T$ on $\AA_\ast(X)$ from
Example~\ref{X-based-duality}; this is no longer considered quite
right for $\AA_\ast(X,w)$. To correct it we define a ``local
coefficient system'' $\Gamma_!$ of infinite cyclic groups on $X$ by
$\Gamma_!(\sigma)=\Gamma\times_{\pi}p^{-1}(\hat\sigma)$, for
simplices $\sigma\in sX$ with barycenter $\hat\sigma$. Then we
compose the old chain duality $T$ with the endofunctor given by
\[   \sum_{\sigma\in sX} M(\sigma) \,\,\mapsto \,\,
\sum_{\sigma\in sX} M(\sigma)\otimes \Gamma_!(\sigma) \] to obtain
the new duality. Hence, for objects $M$ and $N$ of $\AA_\ast(X)$,
there is an embedding
\[ M\otimes_{\AA_\ast(X)}N
\lra C_*(X;\Gamma_!)\otimes(M\otimes_{\ZZ}N). \] The assembly
functor $\alpha \colon \AA_\ast(X,w) \rightarrow \AA[\pi]^w$ is
defined exactly as in the untwisted setting by
\[
M \mapsto \alpha (M) = \sum_{\sigma^\natural \in sX^\natural}
M(\sigma^\natural)
\]
with $M(\sigma^\natural):=M(p(\sigma^\natural))$, but the old
definition of the comparison maps
\[ C\otimes_{\AA_\ast(X)}D \lra \alpha(C)\otimes_{\AA[\pi]}\alpha(D)
\]
for $C,D$ in $\BB(\AA_\ast(X))$ has to be modified since its source
and target are not what they were then. This is straightforward. As
before, the assembly functor induces a map between L-theory spaces
\[
\alpha \colon \bL_n(\AA_\ast (X,w)) \rightarrow \bL_n (\AA[\pi]^w)
\]
also called assembly. Assuming that $X$ is connected and $p\co
X^\natural \to X$ is a universal covering, one can define the space
$\bS(X,n,w)$ as the homotopy fiber of the assembly map. There is
also a description of this space as an $L$-theory space of the
category of chain complexes in $\BB(\AA_\ast(X,w))$ with
contractible assembly in $\BB(\AA[\pi]^w)$.

\

\nin\textbf{Truncated version.} Let $X$ be a finite $\Delta$-complex
with subcomplexes $X_1$ and $X_2$ such that $X_1\cup X_2=X$.
Ranicki's excision theorem for $\bL_{\bullet}(\AA_\ast(X))$
mentioned earlier in Example~\ref{hlgy} means that the square of
$\Omega$--spectra
\[
\xymatrix{
  \bL_{\bullet}(\AA_\ast(X_1\cap X_2)) \ar[d] \ar[r]
&  \bL_{\bullet}(\AA_\ast(X_1)) \ar[d] \\
 \bL_{\bullet}(\AA_\ast(X_2)) \ar[r]
&  \bL_{\bullet}(\AA_\ast(X_1\cup X_2)) }
\]
is a homotopy pushout square. This implies that $X\mapsto
\pi_*\bL_{\bullet}(\AA_\ast(X))$ is a generalized (unreduced)
homology theory. We now recall some of Postnikov's method for making
truncated variants of generalized homology theories.

Let $X\mapsto Q_*(X)$ be a generalized homology theory (from finite
$\Delta$-complexes to graded abelian groups, say). Fix $r\in \ZZ$.
We define the Postnikov fiber truncation $\wp^r$ by
\[  \wp^rQ_n(X) = \im [ Q_n(X^{n-r})\to Q_n(X^{n-r+1}) ]. \]
Then $\wp^rQ_*$ is again a generalized homology theory and there are
long exact sequences
\[ \cdots \to H_{n-r+1}(X;Q_k) \to \wp^{r+1}Q_n(X) \to \wp^rQ_n(X)
\to H_{n-r}(X;Q_r) \to \cdots \] where $Q_r=Q_r(\ast)$. For fixed
$r\le n-\dim(X)$, we clearly have $\wp^rQ_n(X)=Q_n(X)$.

In our situation there is a space version of this construction which
we outline very briefly. Let $(C,\psi)$ be a $k$--simplex in
$\bL_n(\AA_\ast(X))$, alias $n$--dimensional QAPC in
$(\AA_\ast(X))^\ast(\Delta^k)$. We say that $(C,\psi)$ is of type
$\wp^r$ if, for simplices $\sigma\in sX$, $\tau\in \Delta^k$, the
chain complex $C(\sigma,\tau)$ is zero if $|\sigma|-|\tau|>n-r+1$,
and nullbordant as a QAPC of dimension $n-|\sigma|+|\tau|$ in
$\BB(\AA)$ if $|\sigma|-|\tau|=n-r+1$ (with the QAPC structure
determined by $\psi$). The simplices of type $\wp^r$ form a
$\Delta$--subset $\bL_n(\AA_\ast(X);\wp^r)$ of $\bL_n(\AA_\ast(X))$.
Letting $n$ vary, these $\Delta$--subsets can be arranged into a
subspectrum
\[ \bL_{\bullet}(\AA_\ast(X);\wp^r) \subset \bL_{\bullet}(\AA_\ast(X))
  \]
and we have $\bL_{\bullet}(\AA_\ast(X);\wp^r)\simeq X_+\wedge
\wp^r\bL_{\bullet}(\AA)$. The verification, along the lines of
Ranicki's reasoning for $r=-\infty$, is left to the reader.

For us, only the case of $\bL_n(\AA_\ast(X);\wp^1)$ with $\dim(X)=n$
is of interest. In that case $\bL_n(\AA_\ast(X);\wp^1)$ is a
$\Delta$--subset of $\bL_n(\AA_\ast(X))$ determined by a condition
on the $0$--simplices only. A 0--simplex $(C,\psi)$ of
$\bL_n(\AA_\ast(X))$, alias $n$--dimensional QAPC in $\AA_\ast(X)$,
belongs to $\bL_n(\AA_\ast(X);\wp^1)$ if and only if
$(C(\sigma),\psi|\sigma)$ is a nullbordant 0-dimensional QAPC in
$\BB(\AA)$ for every $n$--simplex $\sigma$ in $X$. (The
higher--dimensional simplices in $\bL_n(\AA_\ast(X))$ belong to
$\bL_n(\AA_\ast(X);\wp^1)$ precisely if all their vertices do.) Note
that if $X$ is connected, then the clause \emph{for every
$n$--simplex $\sigma$}  can be replaced by \emph{for some
$n$--simplex $\sigma$}.

\

\nin\textbf{Geometric versus algebraic surgery sequence.} For a
$\Delta$-complex $X$ which is an $n$-dimensional connected oriented
manifold and a specified universal covering of $X$ with deck
transformation group $\pi$, there is the following diagram of
homotopy fibration sequences
\[
\xymatrix{
  \sStw (X) \ar[r] \ar[d] & \sN (X) \ar[r] \ar[d] & \sL_n (X) \ar[d] \\
  \bS(X;\wp^1) \ar[r] & \bL_n(\AA_\ast (X);\wp^1) \ar[r] & \bL_n (\AA[\pi])
}
\]
where the vertical arrows are homotopy equivalences, and
$\bS(X;\wp^1)$ is \emph{defined} as the homotopy fiber of
$\bL_n(\AA_\ast (X);\wp^1) \to \bL_n (\AA[\pi])$.

Note that in Example \ref{frag-quad-con} we essentially described a
map $\sN (X) \rightarrow \bL_n (\AA_\ast (X))$. This factors through
$\bL_n(\AA_\ast (X);\wp^1)$. Indeed for an $n$-dimensional simplex
$\sigma$ in $X$ the dual cell $D(\sigma,X)$ is a point $\hat\sigma$
and for any degree one map $f\co M\to X$ transverse to $\hat\sigma$,
the inverse image $f^{-1}(\hat\sigma)$ is a closed $0$--manifold of
signature $1$. (Hence the signature of $f^{-1}(\hat\sigma)$ minus
the signature of $\hat\sigma$ is 0.) Again there are versions of the
identification for the cases when $X$ is a non-orientable or just
non--oriented manifold. The details are again left to the reader.

\begin{rem} \label{geo-sym-str}
Strictly speaking, in order to define a map $\sN(X) \ra \bL_n
(\AA_\ast (X))$ we should have added CW-structures as in
example~\ref{diss-duality} and ``geometric symmetric structures''
(maps $\nabla^\sharp$ as in example~\ref{sym-con}) on the manifolds
or CW-spaces involved to the geometric data, since choices of these
must be made before the algebraic data can be extracted. However,
these choices are ``contractible'' choices. Adding them or
neglecting them does not change the homotopy type of $\sN(X)$.
\end{rem}

\nin \textbf{Products.} For $(C,\phi)$ an $m$-dimensional SAPC in
$\BB(\AA[\pi])$, and $(C',\phi')$ an $n$-dimensional SAPC in
$\BB(\AA[\pi'])$, we have $(C \otimes C',\phi \otimes \phi')$~, an
$(m+n)$-dimensional SAPC in $\BB(\AA[\pi \times \pi'])$. For
$(C,\phi)$ an $m$-dimensional SAPC in $\BB(\AA[\pi])$, and
$(D,\psi)$ an $n$-dimensional QAPC in $\BB(\AA[\pi'])$, we have
$(C\otimes D,\phi \otimes \psi)$~, an $(m+n)$-dimensional QAPC in
$\BB(\AA[\pi \times \pi'])$. See \cite[section 8]{RaLMS1}.

Let $X$ be an $m$-dimensional, and $Y$ an $n$-dimensional Poincar\'e
duality CW-complexes with chosen orientation classes. Then we have a
natural isomorphism
\begin{equation} \label{sym-prod}
(C_\ast (\widetilde X \times \widetilde Y),\phi (X \times Y)) \cong
(C_\ast \widetilde X \otimes C_\ast \widetilde Y,\phi(X) \otimes
\phi (Y))
\end{equation}
of $(m+n)$-dimensional SAPCs in $\BB(\AA [\pi_1 (X \times Y)])$. See
\cite[section 8]{RaLMS2}. We are assuming that geometric symmetric
structures on $X$ and $Y$ as in example~\ref{sym-con} have been
selected, and use the product geometric symmetric structure on
$X\times Y$. A similar but slightly more complicated statement for
normal maps and quadratic structures is available. We do not
formulate this because we will not need it, thanks to the stated
equivalence of definitions~\ref{defn-otherstruc}
and~\ref{defn-struc}.


\section{Orthogonal Calculus and Products}
\label{oc-and-prod}

The orthogonal calculus \cite{We} is about continuous functors from
a certain category $\sJ$ of real vector spaces to the category of
spaces. For details and definitions, see also \cite{Ma}. Here we
take another look at orthogonal calculus from a ``multiplicative''
viewpoint.

\medskip
Let $\sJ\iso$ be the subcategory of the isomorphisms in $\sJ$. The
objects of $\sJ\iso$ are the finite dimensional real vector spaces
$V,W,\dots $ with inner product, and the space of morphisms from $V$
to $W$ in $\sJ\iso$ is the space of invertible linear isometries
from $V$ to $W$.

\begin{defn}
\label{defn-products}  Let $E$ and $F$ be continuous functors from
$\sJ\iso$ to based spaces. A \emph{multiplication} on $E$ is a
binatural based map $m\co E(V)\wedge E(W) \to E(V\oplus W)$, defined
for $V$ and $W$ in $\sJ\iso$, which satisfies the appropriate
associativity law. A \emph{unit} for the multiplication is a
distinguished element $1\in E(0)$ which is a neutral element for the
multiplication $m$. For $E$ equipped with a multiplication $m$ and a
unit, an \emph{action} of $E$ on $F$ is a binatural based map \[a\co
E(V)\wedge F(W)\to F(V\oplus W),\] again defined for $V$ and $W$ in
$\sJ\iso$, which satisfies the appropriate associativity law
involving $m$ and $a$, and has $1\in E(0)$ acting by identity maps.
\end{defn}

\begin{expl}
\label{expl-terminalaction}  Let $E(V)=\SS^0$ for all $V$, with
$m\co E(V)\wedge E(W)\cong E(V\oplus W)$ for all $V,W$. For a
continuous $F$ from $\sJ\iso$ to based spaces, an action of $E$ on
$F$ amounts to an extension of $F$ from $\sJ\iso$ to $\sJ$.
\end{expl}

\begin{expl}
\label{expl-initialaction}  Let $E(V)=\ast$ for all $V$. For any
continuous $F$ from $\sJ\iso$ to based spaces, there is a unique
action of $E$ on $F$.
\end{expl}

\begin{expl}
\label{expl-freeaction}  Let $E$ be given as in
definition~\ref{defn-products}, with multiplication $m$. Fix $U$ in
$\sJ\iso$. We are going to define an $F$ from $\sJ\iso$ to spaces,
with an action of $E$, in such a way that $F$ is \emph{free} on one
generator $\iota\in F(U)$. We set $F(W)=\ast$ if $\dim(W)<\dim(U)$.
For $W$ with $\dim(W)-\dim(U)=k\ge 0$ we define
\[ F(W)  = \mor(U\oplus\RR^k,W)_+\wedge_{\Or(k)} E(\RR^k). \]
Here ``$\mor$'' refers to a space of morphisms in $\sJ\iso$, and the
$\wedge_{\Or(k)}$ notation means that we are dividing by the
equivalence relation which identifies $(gh,x)$ with $(g,hx)$
whenever $h\in \Or(k)\subset \Or(U\oplus\RR^k)$. It is clear that
$F$ is a functor on $\sJ\iso$. The multiplication $m$ on $E$
determines an action of $E$ on $F$ as follows. For $x\in E(V)$ we
have the left multiplication $m_x\co E(\RR^k)\to E(V\oplus\RR^k)$
and we define the action $a_x\co F(W)\to F(V\oplus W)$ by
\[
\CD \mor(U\oplus\RR^k,W)_+\wedge_{\Or(k)} E(\RR^k) @. \\ @VV
{\textup{incl.} \wedge m_x} V  @. \\ \mor(U\oplus V\oplus
\RR^k,V\oplus W)_+\wedge_{\Or(V\oplus\RR^k)} E(V\oplus \RR^k)
@>\cong>> F(V\oplus W).
\endCD
\]
The generator of $F$ is $\iota=(\id,1)\in F(U)$. \newline Let $F_1$
be another continuous functor from $\sJ\iso$ to spaces with an
action of $E$. Then a map $v\co F\to F_1$ which respects the actions
of $E$ is completely determined by $v(\iota)\in F_1(U)$, which can
be prescribed arbitrarily.
\end{expl}

\begin{defn}
\label{defn-freeaction}  Given $E$ and $F$ as in
definition~\ref{defn-products}, with multiplication $m$ and action
$a$, we say that $F$ is \emph{free} if it has a wedge decomposition
$F\cong \bigvee F_{\lambda}$ where each $F_{\lambda}$ is free on one
generator, as in example~\ref{expl-freeaction}.
\end{defn}

\begin{defn}
\label{defn-CW}  Let $E,F$ be as in definition~\ref{defn-products},
with multiplication $m$ on $E$ and action $a$ of $E$ on $F$. An
\emph{$E$-$CW$-structure} on $F$ is a collection of subfunctors
$F^i\subset F$ for $i=-1,0,1,2,\dots$, subject to a few conditions:
\begin{itemize}
\item $F^{-1}=\ast$ and $F^i\subset F^{i+1}$ for $i\ge -1$~;
\item $F(V)=\bigcup_i F^i(V)$ with the colimit topology, for all $V$~;
\item the action of $E$ on $F$ respects each $F^i$~;
\item for every $i\ge -1$, there exists a pushout square
\[
\CD Z^i\wedge\SS^i_+ @>\subset >> Z^i\wedge\DD^{i+1}_+ \\ @VVV  @VVV
\\ F^i @>\subset >> F^{i+1}
\endCD
\]
where $Z^i$ is another functor from $\sJ\iso$ to spaces, with a free
action of $E$, and the (vertical) arrows respect the actions of $E$.
\end{itemize}
The subfunctor $F^i$ is sometimes called the $i$-skeleton of $F$.
\end{defn}

\begin{lem}
\label{lem-CWapprox} Let $E,F$ be as in
definition~\ref{defn-products}, with multiplication $m$ on $E$ and
action $a$ of $E$ on $F$. There exists an $E$-$CW$-approximation for
$F$. That is, there exists a weak equivalence $\hat F\to F$ of
continuous functors on $\sJ\iso$ with $E$-action, where $\hat F$ has
an $E$-$CW$-structure as in~\ref{defn-CW}.
\end{lem}

\begin{defn}
\label{defn-inducing}
  Let $E_1,E_2$ be functors as in
definition~\ref{defn-products}, with multiplications $m_1$ on $E_1$
and $m_2$ on $E_2$. Let $h\co E_1\to E_2$ be a natural
transformation respecting the units and multiplications. We consider
continuous functors $F$ from $\sJ\iso$ to pointed spaces, either
with an action of $E_1$ or with an action of $E_2$. Composition with
$h$ gives a functor from the category of functors $F$ as above with
an action of $E_2$ to the category of functors $F$ as above with an
action of $E_1$. This functor has a left adjoint, which we call
\emph{induction along $h$} and denote by $\ind_h$. Thus, for $F$
from $\sJ\iso$ to based spaces with an action of $E_1$~, we have
$\ind_hF$ from $\sJ\iso$ to based spaces with an action of $E_2$.
There is a canonical transformation $F\to \ind_hF$ which
``intertwines'' the actions and has a universal property.
\end{defn}

The induction functor $\ind_h$ as defined above tends to produce
pathological results. However, there are situations where it is well
behaved:

\begin{lem}
\label{lem-goodinduction} Keep the notation of
definition~\ref{defn-inducing}. Suppose that $F$ from $\sJ\iso$ to
based spaces comes with an action of $E_1$. If $F$ has an
$E_1$-$CW$-structure with skeletons $F^i$~, then $\ind_hF$ has an
$E_2$-$CW$-structure with skeletons $\ind_hF^i$.
\end{lem}

We therefore have something like a ``derived induction'' procedure
which is as follows. Fix $h\co E_1\to E_2$ as above and some $F$
from $\sJ\iso$ to based spaces, with an action of $E_1$. Replace $F$
by an $E_1$-$CW$-approximation as in lemma~\ref{lem-CWapprox}. Then
apply $\ind_h$ to the $CW$-approximation, assuming that a
multiplicative $h\co E_1\to E_2$ is given. \newline Our interest
here is mainly in the case where $E_2(V)=\SS^0$ for all $V$ as in
example~\ref{expl-terminalaction}. Then $\ind_h$ of an
$E_1$-$CW$-approximation to $F$ is a continuous functor on $\sJ$,
and that is (still) the sort of object we are after.

\begin{lem}
\label{lem-boring} Keep the assumptions of
lemma~\ref{lem-goodinduction}. Suppose in addition that $h\co
E_1(V)\to E_2(V)$ is a based homotopy equivalence for every $V$.
Then the canonical map $F(W)\to \ind_hF(W)$ is a based homotopy
equivalence for every $W$.
\end{lem}

\begin{defn}
\label{joincat} Let $\sP$ be the following monoidal category.
Objects of $\sP$ are pairs $(X,u)$ where $X$ is a finite
$\Delta$-complex homotopy equivalent to a sphere and $u\colon X\to
X$ is a free involution respecting the $\Delta$-complex structure.
Morphisms are $\Delta$-maps, respecting the involutions, which are
weak homotopy equivalences. The monoidal operation $*$ is the join:
$(X,u)*(Y,v)=(X*Y,u*v)$.
\end{defn}

\medskip
\emph{Comment.} If $X$ is a finite Delta-complex homotopy equivalent
to $S^m$ and $u$ is a free involution on $X$, then by the Lefschetz
trace formula $u$ acts on the reduced $m$-th homology of $X$ by
$(-1)^{m+1}$. Using that observation and obstruction theory, one can
easily show that the orbit space $X/u$ is homotopy equivalent to
$\RR P^m$. The case $m=-1$ is \emph{not} an exception: in that case
$X=\emptyset$ and the reduced $(-1)$-th homology is $\cong \ZZ$.
Obviously $\emptyset$ is a very important object of $\sP$ because it
is a unit for the join operation.

\medskip
We think of $\sP$ as a combinatorial variant of the monoidal
category $\sJ\iso$ (the subcategory of $\sJ$ in which only
isomorphisms are allowed as morphisms, with the monoidal operation
\emph{product} alias direct sum). The following definitions
introduce a construction, essentially a homotopy Kan extension,
which ``transforms'' a space-valued functor on $\sP$ into a
space-valued continuous functor on $\sJ\iso$.

\begin{defn} For $V$ in $\sJ\iso$, let $\sP_V$ be the following topological
category. An object is an object $(X,u)$ of $\sP$ together with a
map $\lambda\colon X\to V\smin 0$ which is linear on each simplex of
$X$ and satisfies $\lambda u(x) =-\lambda(x)$. A morphism from
$(X,u)$ with $\lambda_0\colon X\to V$ to $(Y,v)$ with
$\lambda_1\colon Y\to V$ is a morphism $f\colon (X,u)\to (Y,v)$ in
$\sP$ satisfying $\lambda_1f=\lambda_0$. The topology on the object
class $\sP_V$ comes (only) from the fact that, in an object $(X,u)$
with $\lambda\colon X\to V$, the $\lambda$ can vary continuously
(within a finite dimensional space, since $\lambda$ is determined by
its values on vertices of $X$). Thus the projection functor
$\sP_V\to \sP$ is continuous.
\end{defn}

\begin{defn}
\label{defn-Kanext} For a functor $G$ from $\sP$ to (well-)based
spaces, let $G^K$ be the continuous functor on $\sJ\iso$ defined by
\[  G^K(V) = \hocolimsub{\lambda\colon X\to V} G(X,u) \]
where $\lambda\colon X\to V$ runs through $\sP_V$ (and the homotopy
colimit is ``reduced'' so that it is again a based space).
\end{defn}

\begin{lem}
\label{lem-howtosee} Suppose that $G$ takes all morphisms in $\sP$
to homotopy equivalences. Then for any $\lambda\colon X\to V$ in
$\sP_V$ the inclusion $G(X,u)\to G^K(V)$ is a homotopy equivalence.
\end{lem}

\proof The hypothesis implies that the forgetful map from $G^K(V)$
to $B\sP_V$ is a quasi-fibration. Therefore it is enough to show
that $B\sP_V$ is contractible. To do so we first replace the
topological category $\sP_V$ by a simplicial category $k\mapsto
\sP_{V,k}$. An object of $\sP_{V,k}$ is an object $X=(X,u)$ of $\sP$
together with a map
\[ f\co X\times\Delta^k\lra V \]
such that, for every $z\in \Delta^k$, the map $f_z\co X\to V$
defined by $f_z(x)=f(x,z)$ defines a n object in $\sP_V$. A morphism
from $f\co X\times\Delta^k\to V$ to $g\co Y\times\Delta^k\to V$ is a
$\Delta$-map $X\to Y$ making a certain triangle commute. It is easy
to show that, for fixed $j\ge 0$, the canonical map
\[  \left|k\mapsto N_j\sP_{V,k}\right| \lra N_j\sP_V \]
is a homotopy equivalence, where $N_\bullet$ denotes the nerve
construction. Integrating over $j$, we conclude that the canonical
map from $\left|k\mapsto B\sP_{V,k}\right|$ to $\lra B\sP_V$ is a
homotopy equivalence. Now it only remains to show that $B\sP_{V,k}$
is contractible for fixed $k\ge 0$. But that is again
straightforward because $\sP_{V,k}$ is a \emph{directed} category,
in the strong sense that every finitely generated diagram in
$\sP_{V,k}$ admits a co-cone. Namely, suppose that $\sD$ is a
finitely generated category and $v\co \sD\to \sP_{V,k}$ is any
functor. Then there exists a constant functor $c\co \sD\to
\sP_{V,k}$ and a natural transformation $v\Rightarrow c$. To
construct $c$, form the ``direct limit'' of $v$, not necessarily an
object of $\sP_{V,k}$ but a well defined finite $\Delta$-complex $Z$
with a certain map from $Z\times\Delta^k$ to $V$ avoiding $0\in V$.
By attaching simplices to $Z$ as appropriate, embed $Z$ in an object
of $\sP_{V,k}$. \qed

\begin{lem}
\label{lem-howtomultiply} Let $G_0$, $G_1$ and $G_2$ be functors
from $\sP$ to (well-)based spaces. Any natural multiplication
\[ G_0(X,u)\wedge G_1(Y,v)\lra G_2(X*Y,u*v) \]
induces a natural multiplication $G_0^K(V)\wedge G_1^K(W) \lra
G_2^K(V\times W)$.
\end{lem}


\section{Joins in dissected $L$-theory} \label{alg-join}

Let $\AA$ be the additive category of finitely generated free
abelian groups, with the standard chain duality. Let $(C,D)$ and
$(C',D')$ be chain complex pairs in $\BB(\AA)$. We assume the
boundary inclusions $D\to C$ and $D'\to C'$ are cofibrations, i.e.,
degreewise split, and to be quite precise we assume that such
degreewise splittings have been specified. Then there is the product
pair $(C,D)\otimes(C',D')$ consisting of $C\otimes C'$ and the
subcomplex $(C\otimes D')\oplus_{(D\otimes D')}(D\otimes C')$ as the
``boundary''.

Let $X$ and $Y$ be finite $\Delta$-complexes. By a \emph{dissection
of $D$ over $X$} we mean a splitting $D=\sum_{\sigma\in
sX}D(\sigma)$ of $D$ as a graded abelian group which promotes $D$ to
an object of $\BB(\AA_\ast(X))$.

\begin{lem}
\label{lem-fragment} A dissection of $D$ over $X$ and a dissection
of $D'$ over $Y$ together determine a dissection of the boundary of
$(C,D)\otimes(C',D')$ over $X*Y$.
\end{lem}

\proof Suppose that the dissections of $D$ and $D'$ are given by
graded group splittings $D=\bigoplus_\sigma D(\sigma)$ and
$D'=\bigoplus_\tau D(\tau)$. The specified splittings of $D\to C$
and $D'\to C'$ also give us identifications $C\cong D\oplus C/D$ and
$C'\cong D'\oplus C'/D'$ of graded groups. Hence the boundary
complex of $(C,D)\otimes(C',D')$ splits (as a graded group) into
summands
\[ D(\sigma)\otimes D'(\tau)~,\quad D(\sigma)\otimes C'/D'~,\quad C/D\otimes
D'(\tau)\,. \] We now label these summands by simplices of $X*Y$. A
summand of the form $D(\sigma)\otimes D'(\tau)$ gets the label
$\sigma*\tau$. A summand of type $D(\sigma)\otimes C'/D'$ gets the
label $\sigma$, and we note that $X\subset X*Y$. A summand of type
$C/D\otimes D'(\tau)$ gets the label $\tau$. It is easy to verify
that this labeling defines a dissection. \qed

\medskip
We next discuss a few variations on the theme of
lemma~\ref{lem-fragment} where the pairs $(C,D)$ and $(C',D')$ come
with symmetric structures. \emph{Notation}: We write for example
\[ (D\otimes D)^{h\ZZ/2}~\quad, \quad~(D\otimes D)_{h\ZZ/2} \]
for $\Hom_{\ZZ[\ZZ/2]}(W,D\otimes D)$ and
$W\otimes_{\ZZ[\ZZ/2]}(D\otimes D)$, respectively, where $W$ is the
standard $\ZZ[\ZZ/2]$-resolution of $\ZZ$. A symmetric structure on
a pair $(C,D)$ of chain complexes will be described as a chain
$\varphi\in (C\otimes C)^{h\ZZ/2}$ whose boundary $\partial\varphi$
is in the subcomplex $(D\otimes D)^{h\ZZ/2}$. Similarly a quadratic
structure on $(C,D)$ will be described as a chain $\psi\in (C\otimes
C)_{h\ZZ/2}$ whose boundary $\partial\varphi$ is in the subcomplex
$(D\otimes D)_{h\ZZ/2}$.

\begin{lem}
\label{lem-morefragment} Keep the assumptions of
lemma~\ref{lem-fragment}. Suppose also that the pairs $(C,D)$ and
$(C',D')$ are equipped with symmetric structures $\varphi$ and
$\psi$. If the boundary symmetric structures $\,\partial\varphi$ on
$D$ and $\,\partial\psi$ on $D'$ are dissected over $X$ and $Y$,
respectively, then the symmetric structure
$\partial(\varphi\otimes\psi)$ on the boundary chain complex of
$(C,D)\otimes(C',D')$ is dissected over $X*Y$. If in addition
$\partial\varphi$ and $\partial\psi$ are dissected Poincar\'e, then
$\partial(\varphi\otimes\psi)$ is dissected Poincar\'e.
\end{lem}

\begin{lem}
\label{lem-yetmorefragment} Keep the assumptions of
lemma~\ref{lem-morefragment}. Suppose also that the
$\Delta$-complexes $X$ and $Y$ come with free actions of a finite
group $\pi$, and the chain complex pairs $(C,D)$, $(C',D')$ come
with actions of $\pi$ so that the dissections of $D$ and $D'$ are
$\pi$-invariant (in the sense that $gD(\sigma)=D(g\sigma)$ and
$gD'(\tau)=D'(g\tau)$ for $g\in \pi$ and simplices $\sigma$ in $X$,
$\tau$ in $Y$). Suppose further that $\varphi$ and $\psi$ are
$\pi$-invariant. Then $\varphi\otimes\psi$ is $\pi$-invariant, and
the dissection of $\partial(\varphi\otimes\psi)$ over $X*Y$ is
$\pi$-invariant for the diagonal action of $\pi$ on $X*Y$.
\end{lem}

We will need slight generalizations of
lemmas~\ref{lem-fragment},~\ref{lem-morefragment}
and~\ref{lem-yetmorefragment}. Let $\AA_1$, $\AA_2$ and $\AA_3$ be
additive categories with chain duality. We assume given a functor
\[ \begin{array}{ccc}
\BB(\AA_1)\times \BB(\AA_2) & \lra & \BB(\AA_3) \\ (C,C') & \mapsto
& C\boxtimes C'
\end{array}
\]
which is bi-additive and respects cofibration sequences in any of
the two input variables. (This means that if $C$ appears in a
degreewise split short exact sequence $K\to C\to Q$, then
$K\boxtimes C'\to C\boxtimes C'\to Q\boxtimes C'$ is also degreewise
split short exact, and similarly if $C'$ appears in a degreewise
split short exact sequence.) We also need some compatibility between
$\boxtimes$ and the tensor products $\otimes_{\AA_i}$ for $i=1,2,3$.
We assume therefore that a natural chain map
\[  u\co (B\otimes_{\AA_1}C) \otimes (B'\otimes_{\AA_2}C')
\lra (B\boxtimes B')\otimes_{\AA_3}(C\boxtimes C') \] is given,
depending on variables $B,C$ in $\BB(\AA_1)$ and $B',C'$ in
$\BB(\AA_2)$. This is supposed to respect nondegenerate cycles. That
is, if $x$ and $y$ are nondegenerate cycles in $B\otimes_{\AA_1}C$
and $B'\otimes_{\AA_2}C'$~, respectively, then $u(x,y)$ is a
nondegenerate cycle in $(B\boxtimes B')\otimes_{\AA_3}(C\boxtimes
C')$.
\newline
Given $C$ in $\BB(\AA_1)$ and $C'$ in $\BB(\AA_2)$ and symmetric
structures $\varphi$, $\psi$ on $C$ and $C'$~, respectively (of
degrees $m$ and $n$, respectively) we have a symmetric structure
$\varphi\otimes\psi$ on $C\boxtimes C'$, of degree $m+n$, by
composing
\[
\CD W @>\textrm{\quad diagonal\quad}>> W\otimes W\qquad @. \\ @.
@VV\varphi\times\psi V   @. \\ @.
\quad\qquad(C\otimes_{\AA_1}C)\otimes(C'\otimes_{\AA_2}C') @>u>>
(C\boxtimes C')\otimes_{\AA_3}C\boxtimes C').
\endCD
\]
There is a similar construction for pairs $(C,D)$ and $(C',D')$ with
symmetric structures $\varphi$ and $\psi$, respectively.

\begin{lem}
\label{lem-withparameters}
Lemmas~\ref{lem-fragment},~\ref{lem-morefragment}
and~\ref{lem-yetmorefragment} remain valid for pairs $(C,D)$ and
$(C',D')$ in $\BB(\AA_1)$ and $\BB(\AA_2)$~, respectively, in which
case the product pair is to be taken as $(C,D)\boxtimes(C',D')$ in
$\BB(\AA_3)$.
\end{lem}

\begin{expl}
\label{expl-preview}  We use the above ideas on dissection and joins
to produce examples of functors $F$ and $E$ on $\sP$ satisfying the
condition of lemma~\ref{lem-howtosee} and related by multiplications
as in lemma~\ref{lem-howtomultiply}. Fix $X$ in $\sP$. Determine $m$
so that $X\simeq S^{m-1}$. Let $\AA$ be the category of finitely
generated abelian groups, as before. In outline, $E(X)$ will be
defined as the algebraic cobordism space of formally $m$-dimensional
symmetric Poincar\'e pairs $(C,D,\varphi)$ in $\BB(\AA)$ where the
boundary $(D,\partial\varphi)$ is equipped with an involution and a
(symmetric Poincar\'e) dissection over $X$ which respects the
involutions. The definition of $F(X)$ is the same, except for an
additional condition on the symmetric Poincar\'e pairs
$(C,D,\varphi)$, which is that both $C$ and $D$ have to be
contractible as objects of $\BB(\AA)$. A tensor product
construction, where we use the above lemmas on dissection and joins,
gives us maps
\[ E(X)\wedge E(Y) \lra
E(X*Y)~, \quad E(X)\wedge F(Y)\lra F(X*Y). \] These maps have the
usual associativity properties. There is also a unit in
$E(\emptyset)$. For all details, see the next chapter.
\end{expl}

\begin{rem}  Lemma~\ref{lem-howtosee} gives us some
information about the associated functors $F^K$ and $E^K$ on
$\sJ\iso$. We will use that in the next chapter to understand $F^K$
in homotopy theoretic terms. But we will not attempt to describe the
homotopy type of $E^K(V)$ for all or some $V$. Instead we will use
geometric ideas in a later chapter to construct a ``smaller''
functor $E^{K,\eta}$ on $\sJ\iso$, with multiplication as in
definition~\ref{defn-products}, and a multiplicative transformation
$E^{K,\eta}\to E^K$. The point of the smallness is that $E^{K,\eta}$
admits a weak equivalence $h$ to the constant functor $V\mapsto
\SS^0$ with the standard product. We regard $F^K$ as a functor with
an action of $E^{K,\eta}$ and do a derived induction along $h$ to
obtain a functor defined on all of $\sJ$. See
example~\ref{expl-terminalaction} and lemma~\ref{lem-boring}.
\end{rem}


\section{Structure spaces in the algebraic setting} \label{str-spaces}

We begin with a remark on the additive categories $\AA_*(X)$ and
$\AA^*(X)$ with chain duality. They have been defined for any
$\Delta$-complex $X$. However, it should be clear that the
definitions extend to more general cases where the faces of $X$ are
``convex polytopes''. We will only need this extension in the case
of $\AA^*(X)$, and then only when $X$ is a ``multisimplex'', that
is, a product of finitely many standard simplices.

\begin{defn}
\label{defn-Em} Here we give the full definition of $E(X)$ and
$F(X)$ in example~\ref{expl-preview}. Both are geometric
realizations of $m$-fold $\Delta$-sets. Let $k=(k_1,\dots,k_m)$ be a
multi-index. Write $\Delta^k$ for
$\Delta^{k_1}\times\cdots\times\Delta^{k_m}$ and
\[ E(X,k)\]
for the set of $k$-multisimplices of $E(X)$. By definition, an
element of that set is a pair $(C,D)$ in the category
$\BB(\AA^*(\Delta^k))$. There are more data:
\begin{itemize}
\item We ask for an SAPC pair structure $\psi$ of formal
dimension $m+\sum_ik_i$ on the chain complex pair $(C,D)$.
\item A dissection of $(D,\partial\psi)$ over $X$ is also part
of the data, and we want this to be Poincar\'e; hence the dissected
$(D,\partial\psi)$ is an SAPC in
\[ \BB((\AA^*(\Delta^k))_*(X)) \,
\cong\, \BB((\AA_*(X))^*(\Delta^k))\,. \]
\item An involution on $(C,D)$ is required, respecting
the dissection of $D$ (and compatible with the given free involution
on $X$, as far as $D$ is concerned) and respecting $\psi$ up to a
sign $(-1)^m$.
\end{itemize}
Note that $E(X,k)$ is a pointed set: the zero object $(C=D=0)$
serves as the base point. When we form the geometric realization
\[ \left|\,k\mapsto
E(X,k) \,\right| \] we collapse all base point simplices to a single
point; hence the geometric realization is a pointed space. The
product maps
\[ E(X)\wedge E(Y) \lra
E(X*Y) \] are induced by set maps
\[
\CD E(X,k)\wedge E(Y,\ell) @>>> E(X,k\#\ell )
\endCD \]
with $k\#\ell=(k_1,k_2,\dots,k_m,\ell_1,\ell_2,\dots,\ell_n)$, which
in turn are given by a generalized (but obvious) tensor product
construction $\boxtimes$ as in lemma~\ref{lem-withparameters}.
\newline The definition of $F(X)$ is almost identical with
that of $E(X)$ just given. Again it is the geometric realization of
an $m$-fold $\Delta$-set and we write $F(X,k)$ for the set of
$k$-multisimplices. The elements of that set are pairs $(C,D)$
almost exactly as above, but with one added condition, that $C$ and
$D$ be \emph{contractible} as objects of $\BB(\AA^*(\Delta^k))$.
\end{defn}

There are ``first axis'' subspaces of $E(X)$ and $F(X)$, obtained by
allowing only $k$--multisimplices where $k$ has the form
$(k_1,0,0,\dots,0)$.

\begin{lem}
\label{lem-firstaxis} The inclusions of the first axes in $E(X)$ and
$F(X)$~, respectively, are homotopy equivalences.
\end{lem}

\begin{lem} \label{htpy-type-of-F} Let $F^K$ be the functor on $\sJ\iso$
associated to $F$ as in example~\ref{defn-Kanext}. For $V$ in
$\sJ\iso$ with $\dim(V)=j+1$ there is a homotopy equivalence
\[ F^K(V)\simeq \bS(\RR P(V),j). \]
\end{lem}

\proof The functor $F$ satisfies the condition of
lemma~\ref{lem-howtosee}. Hence $F^K(V)\simeq F(X)$, assuming that
$X$ and $V$ are related as in that lemma. Now use
lemma~\ref{lem-firstaxis} to complete the proof. \qed


\section{Periodicity and Thom isomorphism} \label{proofs}

\begin{thm} We have $F^K(V)\cong \Omega^4F^K(V)$ by a natural
homeomorphism which respects the action of $E^K$.
\end{thm}

\proof The proof is by inspection, using the well-known periodicity
of algebraic $L$-theory given by the double (skew-)suspension. But
it is appropriate to say what exactly $\Omega^4F^K(V)$ means. First
of all, given an $m$-fold based $\Delta$-set $Y$, what do we mean by
$\Omega Y$~? We can define $\Omega Y$ as the based $m$-fold
$\Delta$-set given by
\[ (k_1,k_2,\dots,k_m)
\mapsto \{y\in Y(k_1+1,k_2,\dots,k_m) \mid u^*_{k_1}(y)=\ast,\,\,
v^*_{k_1}(y)=\ast\} \] where $u_{k_1}\co
\{0\}\to\{0,1,\dots,k_1,k_1+1\}$ is the inclusion and
\[ v_{k_1}\co \{0,1,\dots,k_1\}\to\{0,1,\dots,k_1,k_1+1\} \]
is given by $i\mapsto i+1$. (These monotone maps act as face
operators in the ``first'' coordinate direction.) This definition is
justifiable if $Y$ has the Kan extension property.

We defined $F^K(V)$ for $m$-dimensional $V$ as a reduced homotopy
colimit of spaces $F(X)$ for $X\to V$ in $\sP_V$. Since $F$
satisfies the condition of lemma~\ref{lem-howtosee}, we may define
$\Omega^4F^K(V)$ as a reduced homotopy colimit
\[  \hocolimsub{X\to V}\, \Omega^4 F(X)  \]
for $X\to V$ in $\sP_V$. Moreover in that last expression we may
interpret $\Omega^4$ as the fourth power of an operator $\Omega$ on
based $\Delta$-sets, as defined above, to be applied before
geometric realization. Then we have indeed $F^K(V)\cong
\Omega^4F^K(V)$. \qed

\bigskip
Next we discuss Thom isomorphisms. Our theme is that for an oriented
$W$ in $\sJ\iso$ of even dimension $d$, the functor taking $V$ in
$\sJ\iso$ to $\Omega^dF^K(V)$ is identified with a subfunctor of
$V\mapsto F^K(V\oplus W)$, and that the homotopy fibers of this
inclusion are in some sense independent of $V$. Despite appearances,
this has little to do with periodicity.

\medskip
For $X$ in $\sP$ we have a definition of $\Omega^d F(X)$ which looks
almost exactly like the definition of $F(X)$, except for one change
which consists in increasing the formal dimensions of all SAPC's in
the definition of $F(X)$ by $d$. (If $X\simeq S^{m-1}$, then the
$(0,0,...,0)$-simplices of $\Omega^d F(X)$ are certain SAPC pairs
$(C,D,\varphi)$ of formal dimension $m+d$ with a dissection of the
boundary $(D,\partial\phi)$ over $X$.) Using that description of
$\Omega^d F(X)$, we have the following maps:
\[
\xymatrix{
{\Omega^d\big(\!\hocolimsub{f\co X\to V} F(X)\big)} & \\
{\hocolimsub{\twosub{f\co X\to V}{g\co Y\to W}} \Omega^dF(X)\qquad}
\ar[u]_\simeq \ar@<10pt>[r]^\iota & {\hocolimsub{\twosub{e\co Z\to
V\oplus W}{\rule{0mm}{2mm}}} F(Z)}. }
\]
The vertical map is obtained by forgetting the data $g\co Y\to W$,
which run through $\sP_W$, and using the inclusion $\hocolim
\,\Omega^dF\to \Omega^d(\hocolim\, F)$, where the homotopy colimits
are taken over $\sP_V$ only. It is a homotopy equivalence because
its target can be identified up to homotopy equivalence with
$\Omega^dF(X)$ for any $f\co X\to V$ in $\sP_V$ by
lemma~\ref{lem-howtosee}, and its source can also be identified up
to homotopy equivalence with $\Omega^dF(X)$ by the same kind of
argument. The horizontal map uses an embedding
\[ \Omega^dF(X) \lra F(X*Y) \]
for $f\co X\to V$ in $P_V$ and $g\co Y\to W$ in $\sP_W$. This is
simply induced by the inclusions $X\to X*Y$ and $\AA_*(X)\to
\AA_*(X*Y)$, which respect the chain dualities. Using the vertical
arrow as an ``identification'', we write
\[  \iota\co \Omega^dF^K(V) \lra F^K(V\oplus W) \,.\]
We now wish to extend $\iota$ to a homotopy fiber sequence
\[
\CD \Omega^dF^K(V) @>\iota>> F^K(V\oplus W) @>\zeta>> \Phi^K(V\oplus
W,V)
\endCD
\]
where $\Phi^K$ is a certain functor of pairs. The definition of
$\Phi^K$ follows the standard pattern. Hence we start by introducing
functors on pairs of certain $\Delta$-complexes.

\begin{defn}
\label{defn-relativeF} Fix a pair $(Z,X)$ of finite
$\Delta$-complexes, with a free involution, homotopy equivalent to
the pair $(S(V\oplus W),S(V))$ with the antipodal involution. We
define
\[  \Phi(Z,X) \]
essentially by repeating the definition of $F(Z)$ in terms of
Poincar\'e pairs $(C,D,\psi)$ contractible in $B(\AA)$ with
dissected boundary, but relaxing it in one respect: for the
dissection of $\partial\psi$ as a symmetric structure on the
dissected $D$, we require that to be Poincar\'e modulo $X$ only.
(This means that the mapping cone of the appropriate duality map is
chain equivalent, as an object dissected over $Z$, to something
dissected over $X$.) That being done, we put
\[ \Phi^K(V\oplus W,V) := \hocolimsub{(Z,X)\to(V\oplus W,V)}\,
\Phi(Z,X)
\]
The homotopy colimit is taken over all pairs $(Z,X)$ as above and
simplexwise affine $\ZZ/2$-maps $(Z,X)\to(V\oplus W,V)$ taking $X$
to $V\smin 0$ and $Z$ to $V\oplus W\smin 0$, and inducing a homotopy
equivalence $(Z,X)\to(V\oplus W\smin 0, V\smin 0)$. We have maps
\[
\xymatrix{ {\hocolimsub{(Z,X)\to(V\oplus W,V)}\, F(Z)\qquad}
\ar@<5pt>[r]^\zeta \ar[d]_\simeq &
{\hocolimsub{(Z,X)\to(V\oplus W,V)}\, \Phi(Z,X)} \\
{\hocolimsub{Z\to V\oplus W}\, F(Z)} & }
\]
where the horizontal arrow is determined by the inclusions $F(Z)\to
\Phi(Z,X)$. Using the vertical arrow as an identification, we may
write
\[ \zeta\co F^K(V\oplus W) \lra \Phi^K(V\oplus W,W)\,. \]
\end{defn}

\begin{prop}
\label{prop-oldA} For $V$ of dimension $\ge 3$ and oriented $W$ of
even dimension $d$, both in $\sJ\iso$, the following is a homotopy
fiber sequence:
\[
\CD \Omega^dF^K(V) @>\iota>> F^K(V\oplus W) @>\zeta>> \Phi^K(V\oplus
W,V).
\endCD
\]
\end{prop}

\medskip\nin
\emph{Remark.} We need $\dim(V)\ge 3$ to ensure that the inclusion
$V\to V\oplus W$ induces an isomorphism of fundamental groups of the
associated projective spaces.

\medskip
\proof For fixed $V$ and $W$, after routine reformulations as in
lemma~\ref{lem-howtosee}, this boils down to a homotopy fiber
sequence
\[
\CD \Omega^dF(X) @>>> F(X*Y) @>>> \Phi(X*Y,X)
\endCD
\]
where $X\to V$ in $\sP_V$ and $Y\to W$ in $\sP_W$. Both arrows are
inclusion maps. This homotopy fiber sequence is very standard, e.g.
from \cite{Ra}. \qed

\medskip\nin
\emph{Remark.} In the applications of this proposition we need a
fair amount of naturality in the variable $V$. This calls for a more
precise formulation and a better proof. We have the following
commutative square:
\[
\xymatrix{ {\hocolimsub{\twosub{X\to V}{Y\to W}}\,
\Omega^dF(X)\qquad} \ar@<10pt>[r] \ar[d] &
{\hocolimsub{\twosub{(Z,X)\to(V\oplus W,V)}{\rule{0mm}{3mm}}}\, F(Z)} \ar[d] \\
{\hocolimsub{\twosub{X\to V}{Y\to W}}\, \Omega^d\Phi(X,X)\qquad}
\ar@<10pt>[r] & {\hocolimsub{\twosub{(Z,X)\to(V\oplus
W,V)}{\rule{0mm}{3mm}}}\, \Phi(Z,X).} }
\]
Here the upper left-hand term can be identified with
$\Omega^dF^K(V)$ by a forgetful homotopy equivalence which is
natural in $V$ as an object of $\sJ\iso$. The upper right-hand term
can be identified with $F^K(V\oplus W)$, again by a forgetful
homotopy equivalence which is natural in $V$ as an object of
$\sJ\iso$. The lower right-hand term is $\Phi^K(V\oplus W,V)$. The
lower left-hand term is contractible and is functorial in $V$ as an
object of $\sJ\iso$. The square as a whole is homotopy cartesian (by
the argument already given).

\bigskip
Let $E^K$ be the functor on $\sJ\iso$ associated with $E$ of
definition~\ref{defn-Em}, as in definition~\ref{defn-Kanext}. Then
$E^K$ acts on each of the four functors (of the variable $V$)
represented by the four terms of the square just above. For example,
there is an action map
\[ E^K(U) \wedge \Phi^K(V\oplus W,V)\lra \Phi^K(U\oplus V\oplus W,U\oplus V). \]
It is given by a straightforward tensor product construction and we
omit the details.

\medskip
A connected component of $E^K(U)$, with $\dim(U)=m$, determines
(forgetfully) an element in the relative $L^m$ group of the assembly
functor $\AA_*(X) \to \AA_*$~, for any $X\to U$ in $\sP_U$. This
relative $L^m$ group is canonically isomorphic to $L^0(\AA)\cong
\ZZ$. In this way, connected components of $E(U)$ have a ``degree''
which is an integer.

\begin{thm}
\label{thm-oldB} Keep the assumptions of proposition~\ref{prop-oldA}
and let $U$ be in $\sJ\iso$. Let $z\in E^K(U)$ be in a component of
degree $1$. Then multiplication by $z$ is a homotopy equivalence
\[ \Phi^K(V\oplus W,V) \lra \Phi^K(U\oplus V\oplus W,U\oplus V). \]
\end{thm}

\proof Choose $(Z,X)\to (V\oplus W,W)$ as in the definition of
$\Phi^K(V\oplus W,V)$ and choose $Y\to U$ in $\sP_U$. For
typographic reasons we use the abbreviations $Z^Y=Z*Y$ and
$X^Y=X*Y$, and denote passage to $\ZZ/2$-orbits by a tilde
subscript, as in $X_\sim$. We can assume that $z\in E(Y)$, and we
have to show that multiplication by $z$ is a homotopy equivalence
\[ \Phi(Z,X) \lra
\Phi(Z^Y,X^Y). \] Let $k=\dim(V)+d-1=\dim(V)+\dim(W)-1$ and write
$m=\dim(U)$ as before. Most of the proof is in the following
commutative diagram:
\[
\xymatrix{
{\Phi(Z,X)}  \ar[r] &   {\Phi(Z^Y,X^Y)} \\
{\bS(Z_\sim,X_\sim,k)} \ar[u]^\simeq  \ar[r] \ar[d]_\simeq &
{\bS(Z^Y_\sim,X^Y_\sim,k+m)} \ar[u]^\simeq \ar[d]_>>>>>\simeq \\
{\bL_k(\AA_*(X_\sim)\to \AA_*(Z_\sim))} \ar[r] & {\bL_{k+m}\left(\CD
\AA_*(Y_\sim) @>>=> \AA_*(Y_\sim) \\ @VVV @VVV \\ \AA_*(X^Y_\sim)
@>>> \AA_*(Z^Y_\sim)
\endCD\right)}
}
\]
In this diagram , all horizontal arrows are defined as
multiplication by $z\in E(Y)$. The vertical arrows in the upper half
of the diagram are forgetful: they forget boundaries in dissected
Poincar\'e pairs. By proposition~\ref{defn-otherstruc} and
\cite{Ra}, they are homotopy equivalences. That is, the forgotten
dissected boundaries can always be recovered as ``obstructions to
nondegeneracy'' in $\AA_*(Z_\sim)$. To produce the vertical arrow in
the lower half of the diagram, we switch to
definition~\ref{defn-struc} of the algebraic structure spaces, i.e.,
to quadratic structures. (Strictly speaking, we should insert
another row into the diagram to do that.) These vertical arrows in
the lower half of the diagram can then be defined as inclusion maps.
They are homotopy equivalences by the alternative definition of the
algebraic structure spaces as homotopy fibers of assembly maps. Here
we are also exploiting the fact that the inclusions $X_\sim\to
Z_\sim$ and $X^Y_\sim\to Z^Y_\sim$ induce isomorphisms of
fundamental groups, i.e., we are using $\dim(V)\ge 3$. \newline It
remains to show that the lower horizontal arrow in the diagram is a
homotopy equivalence. By a five lemma argument, this reduces to
showing that multiplication by $z$ induces isomorphisms
\[  \pi_*\bL_j(\AA_*(X_\sim)) \quad\lra \quad \pi_*\bL_{j+m}\left(\AA_*(Y_\sim)\to
\AA_*(X^Y_\sim)\right) \] for all $j\in\ZZ$, and similarly with $X$
replaced by $Z$. But this is a case of an ordinary Thom isomorphism.
The standard proof uses a spectral sequence comparison argument. The
spectral sequences are determined by the skeleton filtration of $X$.
\qed

\bigskip
An element of $F^K(V)$, with $\dim(V)=n$, determines forgetfully an
element in the relative $L_n$ group of the assembly functor
$\AA_*(S(V)) \to \AA_*$. This relative $L_n$ group is canonically
isomorphic to $L_0(\AA)$ which we identify with $\ZZ$ using the
isomorphism $\sigma/8$, signature divided by $8$. In this way there
is a degree function from $F^K(V)$ to $\ZZ$ which we denote by
$\tilde\sigma/8$.

\begin{lem}
\label{lem-oldC} Let $U$ be in $\sJ\iso$. Let $z\in E^K(U)$ be in a
component of degree $1$. Then for every $V$ in $\sJ\iso$~, the
following is commutative:
\[
\CD F^K(V) @>z\cdot >> F^K(U\oplus V) \\ @VV {\tilde\sigma/8} V @VV
{\tilde\sigma/8} V  \\ \ZZ  @>=>> \ZZ\,.
\endCD
\]
\end{lem}

\proof The square can be enlarged to a six-term diagram
\[
\CD F^K(V) @>z\cdot >> F^K(U\oplus V) \\ @VV\zeta V       @VV\zeta V \\
\Phi^K(V,0) @>z\cdot >> \Phi^K(U\oplus V,U) \\ @VV {\tilde\sigma/8}
V @VV{\tilde\sigma/8}V  \\ \ZZ  @>=>> \ZZ\,.
\endCD
\]
Here the top square commutes by construction. The middle row is an
isomorphism by theorem~\ref{thm-oldB}, and the reasoning used in the
proof of that theorem also shows that the lower square commutes.
\qed


\section{Structure spaces in the mixed algebraic-geometric setting}
\label{geo-alg-join} In
sections~\ref{alg-sur},~\ref{oc-and-prod},~\ref{alg-join}
and~\ref{str-spaces}, we constructed an algebraic analogue of the
geometric theory in \cite{Ma}, as far as possible.
Section~\ref{proofs} has, in the algebraic setting, the proofs and
theorems that we need in the geometric setting of
section~\ref{proof} and \cite{Ma}. It remains to make the
translation. This is a tedious business and the method that we have
chosen might not be the best. In any case we have in a few places
sacrificed completeness for the sake of intelligibility.

\subsection{Transversality.}
We will rely mainly on the notion of ``transversality to a
foliation''. Let $N$ be a topological space. There is a presheaf on
$N$ which to an open subset $W\subset N$ associates the set of
equivalence relations on $W$. A global section $\rho$ of the
associated sheaf is called a \emph{local equivalence relation} on
$N$. If $N$ can be covered by open subsets $W_{\alpha}$ which admit
homeomorphisms $(p_\alpha,q_\alpha)\co W_\alpha\to V_\alpha\times
U_\alpha$ with $U_\alpha$ open in $\RR^k$ (but no conditions on
$V_{\alpha}$ other than ``being a space'') such that $\rho|W_\alpha$
is represented by the equivalence relation
\[ y\sim z\,\,\,\,\Leftrightarrow\,\,\,\, q_\alpha(y)=q_\alpha(z) \]
on $W_\alpha$~, then $\rho$ is a \emph{codimension $k$ foliation} of
$N$. See \cite{KM} for more details. A map $f$ from a topological
manifold $M$ to $N$ is \emph{transverse} to a codimension $k$
foliation $\rho$ on $N$ if for every $x\in M$ there exists an open
neighborhood $W$ of $f(x)\in N$ and a product structure $(p,q)\co
W\to V\times U$ with $U\subset \RR^k$ representing $\rho|W$, as
above, such that $qf$ is a topological submersion $f^{-1}(W)\to U$.
(To preclude misunderstandings, we point out that submersions don't
have to be surjective. A map between topological manifolds is a
submersion if it satisfies a certain regularity condition at or near
every point of the source manifold.)

Closely related is the following concept of transversality: Suppose
that $N$ is a space, $X\subset N $ is a locally closed subspace, $U$
is an open neighborhood of $X$ in $N$ which comes with a codimension
$k$ foliation, and $X$ is a leaf of that foliation. We say
informally that a map $f$ from a manifold $M$ to $N$ is transverse
to $X$ if $f$ restricted to a sufficiently small neighborhood of
$f^{-1}(X)$ is transverse to the foliation of $U$. In the examples
that we will be looking at, the foliation of $U$ is determined by a
single map $q\co U\to \RR^k$, so that the leaves of the foliation
are the fibers of $q$.

\subsection{Main examples}

\begin{defn} For a $\Delta$-complex $X$, the open cone $\sO(X)$ is
$X_+\wedge[0,\infty)$. We describe a point in $\sO(X)$ as $x=ty$
where $y \in sX$ and $t \in [0,\infty)$, or by its barycentric
coordinates,
\[
x=(x_i)_{i \in \sigma}=(ty_i)_{i\in\sigma}
\]
where $x_i \geq 0$ and $i$ runs through the vertices of a simplex
$\sigma$ containing $y$. The open cone $\sO (X)$ comes with the
\emph{norm} function $x \mapsto \max_i\{x_i\}$. The levels of the
norm define a foliation on $\sO(X)\smin 0$ with the leaves $\sO(X,c)
= \{x \in \sO(X) \; | \; \|x\| = c\}$ for $c>0$. For $(X,u)$ in
$\sP$, the open cone $\sO(X)$ also comes with an involution
$ty\mapsto t \cdot u(y)$, where $y\in sX$.
\end{defn}

\begin{defn}
\label{stratification} We often regard $\sO(X)$ as a stratified
space, stratified by the interiors of the coned dual cells and the
cone point. This is similar to the stratification of $X$ in
example~\ref{diss-duality}. But we need a few more details here. We
need to declare what it means for a map $f\co M\to \sO(X)\smin 0$ to
be transverse to the stratification.

For a $k$-simplex $\sigma$ of $X$ the stratum $\sO(X,\sigma)$ (the
interior of the coned dual cell) consists of points $x\in \sO (X)$
whose barycentric coordinates satisfy $x_i = \|x\|$ if $i \in
\sigma$ and $x_i < \|x\|$ if $i \notin \sigma$. The map $q_\sigma$
defined on a (sufficiently small) neighborhood of $\sO(X,\sigma)$ in
$\sO(X)$ by
\[  x \mapsto \left(\frac{x_i}{\|x\|}\right)_{i\in\sigma}\in \RR^{|\sigma|+1} \]
has image contained in the hypersurface
\[ Z_k=\big\{\,(y_i)_{i\in\sigma}\in \RR^{|\sigma|+1}\mid \max_i\,\{y_i\}=1\}. \]
The sets $q_{\sigma}^{-1}(y)$ for $y\in Z_k$ are the leaves of a
foliation of the neighborhood. One of these leaves is
$\sO(X,\sigma)$ itself. We say that a map $f$ from a manifold to
$\sO(X)$ is \emph{transverse} to $\sO(X,\sigma)$ if $q_{\sigma}f$~,
as a map with target $Z_k$~, is a topological submersion in a
neighborhood of $f^{-1}(\sO(X,\sigma))$.
\end{defn}

\emph{Comment.} The codimension $k$ foliation defined by $q_\sigma$
(of an open neighborhood of $\sO(X,\sigma)$ in $\sO(X)$) also
restricts to a codimension $k$ foliation of an open neighborhood of
$\sO(X,\sigma)\cap \sO(X,c)$ in  $\sO(X,c)$, for every $c>0$.

\medskip
\begin{defn} For $(X,u)$ in $\sP$ with $X\simeq S^{m-1}$ we define $G_0(X,u)$
as an $m$-fold simplicial set. A \emph{nontrivial} $(0,0,\dots,0)$-simplex consists of
a based space $W$ homeomorphic to $\RR^m$, with an involution fixing
the base point, and an equivariant based proper map $p\colon
W\to\sO(X)$ of degree $\pm1$ which is transverse to the foliation of
$\sO(X)\smin 0$ by norm levels, and has preimage of base point equal
to base point. There is also a unique trivial $(0,0,\dots,0)$-simplex which,
along with all its degeneracies, constitutes a connected component of $G_0(X,u)$
after realization.
\end{defn}

\begin{defn} For $(X,u)$ in $\sP$ with $X\simeq S^{m-1}$ we define
$G_1(X,u)$ as an $m$-fold simplicial set. A $(0,0,\dots,0)$-simplex
consists of two based spaces $W$, $W'$ both homeomorphic to $\RR^m$,
both with an involution fixing the base point, and equivariant
proper maps $W\to W' \to \sO(X)$ of degree $\pm1$ such that both
$W'\to \sO(X)$ and the composite map $W\to\sO(X)$ are transverse to
the foliation of $\sO (X)\smin 0$ by norm levels, and have preimage
of base point equal to base point. We regard such a $(0,0,\dots,0)$-simplex as
\emph{trivial} if $W\to W'$ is a homeomorphism. All trivial $(0,0,\dots,0)$-simplices
are to be identified with each other.
\end{defn}

With the above definitions of $G_0$ and $G_1$, there are products
\[
\begin{array}{ccc}
G_0(X,u)\wedge G_0(Y,v)& \lra & G_0(X*Y,u*v)~, \\
G_0(X,u)\wedge G_1(Y,v)& \lra & G_1(X*Y,u*v).
\end{array}
\]
In more detail, on non-trivial simplices the multiplication map is
given simply by the product map, the action map sends $V
\xrightarrow{p} \sO(X)$ and $W \xrightarrow{q} W' \xrightarrow{p'}
\sO(Y)$ to
\[
\xymatrix{ V \times W \ar[r]^{\id \times q} & V \times W' \ar[r]^{p
\times p'} & \sO(X \ast Y). }
\]
They induce similar products involving $G_0^K$ and $G_1^K$.

\medskip
Next we introduce certain refinements of $G_0$ containing more
combinatorial information. These refinements make up a diagram of
functors on $\sP$ and natural transformations
\[  G_{0-3}\to
G_{0-2} \to G_{0-1} \hookrightarrow G_0
\]
Moving from right to left, we first add a transversality condition,
then certain CW-approximations, then cellular diagonal
approximations and cellular fundamental cycles/chains for the
approximating CW-spaces involved. Fix $(X,u)$ in $\sP$, where
$X\simeq S^{m-1}$, and a nontrivial $(0,0,\dots,0)$-simplex $p\colon W\to
\sO(X)$ of $G_0(X,u)$.

\begin{defn}
\label{defn-G_01} To promote $p\colon W\to \sO(X)$ to a
$(0,0,\dots,0)$-simplex in the $m$-fold $\Delta$-set $G_{0-1}(X,u)$,
we impose the condition that $p$ be \emph{simultaneously transverse}
to the strata $\sO(X,\sigma)$ and to the levels of the norm
fibration (see the comment just below). This implies that for every
$c>0$, the restriction of $p$ to the sphere $p^{-1}(c)$ is a map to
$\sO(X,c)\cong X$ which is transverse to the strata $X(\sigma)$.
\end{defn}

\emph{Comment.} Let $\sigma$ be a $k$-simplex of $X$. We test $p~$
for transversality to the stratum $\sO(X,\sigma)$ by asking whether
$q_{\sigma}p$ is a submersion (see definition~\ref{stratification}).
We test for transversality to the norm levels by asking whether
$\|p\|\co W\smin p^{-1}(0)\to \RR$ is a submersion. Here we need a
condition which is slightly stronger than these two transversality
conditions put together. We require that the map
\[ \textup{neighborhood of }p^{-1}(\sO(X,\sigma))
\lra Z_k\times \RR \] defined by $w\mapsto
(q_{\sigma}p(w),\|p(w)\|)$ be a submersion. Since $Z_k\times\RR$ can
be identified with $\RR^{|\sigma|+1}$~, the formula $w\mapsto
(q_{\sigma}p(w),\|p(w)\|)$ can also be replaced by the much simpler
formula
\[   w\mapsto \big(p_i(w)\big)_{i\in\sigma} \in \RR^{|\sigma|+1} \]
where the $p_i(w)$ are the barycentric coordinates of $p(w)$
corresponding to the vertices $i$ of $\sigma$. \newline We mention
the following in passing. Suppose that $p$ satisfies the above
strong transversality condition for a particular $\sigma$. Let
$\tau$ be a face of $\sigma$. Then, in a sufficiently small
neighborhood of $p^{-1}(\sO(X,\sigma))$, the strong transversality
condition for $p$ in relation to $\tau$ and the stratum
$\sO(X,\tau)$ is automatically satisfied. The reason is, of course,
that the barycentric coordinates $p_i(w)$ for $w$ in $W$ and $i$ a
vertex of $\tau$ are subsumed in the barycentric coordinates
$p_i(w)$ for $i$ a vertex of $\sigma$.

\begin{lem}
\label{funproduct} The product on $G_0$ can be refined to a product
on $G_{0-1}$.
\end{lem}

\proof Given $(X,u)$ and $(Y,v)$ in $\sP$ and $p\colon W\to \sO(X)$
and $q\colon W'\to Y$ satisfying the appropriate transversality
conditions, we verify that the resulting map from $W\times W'$ to
$\sO(X)\times\sO(Y)\cong \sO(X*Y)$ satisfies the appropriate
transversality condition. The strata of $\sO(X*Y)\cong
\sO(X)\times\sO(Y)$, apart from the cone point, can be described as
follows:
\begin{itemize}
\item[(i)] for each $\sigma$ in $X$, there is a stratum consisting of all $(x,y)$ in
$\sO(X,\sigma)\times\sO(Y)$ where $\|x\|>\|y\|$~;
\item[(ii)] for each $\tau$ in $Y$, there is a stratum consisting of all $(x,y)$
in $\sO(X)\times\sO(Y,\tau)$ where $\|x\|<\|y\|$~;
\item[(iii)] for each simplex of the form $\sigma*\tau$ with $\sigma$ in $X$ and
$\tau$ in $Y$, there is a stratum consisting of all $(x,y)$ in
$\sO(X,\sigma)\times\sO(Y,\tau)$ where $\|x\|=\|y\|$.
\end{itemize}
Hence the transversality properties that we need follow from the
transversality properties of $p$ in case (i), from the
transversality properties of $q$ in case (ii), and from the
transversality properties of both $p$ and $q$ in case (iii). We omit
the details, except for pointing out that in the case (iii), the
vertex set of $\sigma*\tau$ is identified with the disjoint union of
the vertex sets of $\sigma$ and $\tau$ respectively. It follows that
an expression such as
\[ \big(\big(p\times q)_i(w,w')\big)_{i\in \sigma*\tau} \in \RR^{|\sigma*\tau|+1}\]
for $(w,w')\in W\times W'$ can be re-arranged to look like
\[
\qquad\qquad\big(\big(p_i(w)\big)_{i\in\sigma},\big(q_i(w')\big)_{i\in\tau}\big)\in
\RR^{|\sigma|+1}\times \RR^{|\tau|+1}\,.\qquad\qquad\qed
\]

\medskip
We enlarge $\cat(X)$, the category of simplices of $X$, to a
category $\{0\}*\cat(X)$ by adding the object $0$, its identity
morphism, and one morphism $0 \to \sigma$ for each $\sigma\in sX$.
The transversality condition in the previous definition yields, for
every $(0,0,\dots,0)$-simplex $p\colon W\to \sO(X)$ in
$G_{0-1}(X,u)$ as above, a contravariant functor $W_\diamond$ from
$\{0\}*\cat(X)$ to compact spaces by
\[
\begin{cases}
W_\diamond(\sigma)= W[\sigma][1] \\
W_\diamond(0)= W[0,1]
\end{cases}
\]
where $W[\sigma][1]$ is $p^*$ of $\sO (X,\sigma) \cap \sO (X,1)$,
the norm level $1$ of the coned dual cell corresponding to $\sigma$,
and $W[0,1]$ denotes $p^{-1} \sO (X,[0,1])$, the inverse image of
the portion of $\sO (X)$ with the norm $\leq 1$. (Here $p^*$ denotes
a pullback. We do not write $p^{-1}$ since the dual cell
corresponding to $\sigma$ need not be a subspace of $X$.)

\begin{defn}
\label{G0transCW} To promote $p\colon W\to \sO(X)$ further to a
$(0,\dots,0)$-simplex in $G_{0-2}(X,u)$ we add the following: a
contravariant CW-functor $W_\ddiamond$ from $\{0\}*\cat(X)$ to
compact spaces, with a natural transformation $\gamma\colon
W_\ddiamond\to W_\diamond$ which evaluates to a homotopy equivalence
for every object of $\{0\}*\cat(X)$.
\end{defn}

The definition gives us in particular a $CW$-pair
$(W_\ddiamond(0),W_\ddiamond[1])$ with dissected boundary, where
\[ W_\ddiamond[1]=\colimsub{\sigma\ne 0}W_\ddiamond(\sigma). \]
Passage to cellular chain complexes transforms
a pair of chain complexes $(C,D)$ with $C=C_*(W_\ddiamond(0))$ in
$\BB(\AA)$ and with dissected boundary $D=C_*(W_{\ddiamond}[1])$ in
$\BB(\AA_*(X))$.

\begin{defn}
\label{G0transCWandall} To promote $p\colon W\to \sO(X)$ further to
a $(0,0,\dots,0)$-simplex in $G_{0-3}(X,u)$, we add the following
data: cellular diagonal approximations and fundamental cycles/chains
in the cellular chain complex(es) of $W_\ddiamond$.
\end{defn}

The additional data in definition~\ref{G0transCWandall} imply a
preferred structure of an $n$-dimensional SAPP (where $n = \dim
(W)$) on the pair of chain complexes
\[ (C_*(W_\ddiamond(0)),C_*(W_{\ddiamond}[1])) \]
in $\BB(\AA)$, refined to a dissected $(n-1)$-dimensional SAPC
structure on the boundary in $\BB(\AA_*(X))$.

\medskip
All the above refinements of $G_0$ have multiplications refining the
one on $G_0$. The case of $G_{0-1}$ has already been discussed. For
the case of $G_{0-2}$, suppose given $(X,u)$ and $(Y,v)$ in $\sP$,
as well as $p\colon W\to \sO(X)$ and $q\colon W'\to \sO(Y)$ and
\[ W_\ddiamond\to W_\diamond~, \qquad W'_\ddiamond\to W'_\diamond \]
satisfying the conditions of definition~\ref{G0transCW}. We need to
say what $(W\times W')_\ddiamond$ should be. We have $\{0\} \ast
\cat(X*Y) \cong \{0\} \ast \cat (X) \times \{0\} \ast \cat(Y)$.
Using this identification we let
$(W\times W')_\ddiamond(i,j)=W_\ddiamond(i)\times
W'_\ddiamond(j)$
for $i$ in $\{0\} \ast \cat (X)$ and $j$ in $\{0\} \ast \cat (Y)$.

\begin{lem} The forgetful maps
\[  G_{0-3}\to
G_{0-2} \to G_{0-1} \hookrightarrow G_0
\]
are (weak) homotopy equivalences.
\end{lem}

\medskip
Next, there are refinements of $G_1$ analogous to the above
refinements of $G_0$. These make up a diagram of functors and
natural transformations
\[  G_{1-3}\to
G_{1-2} \to G_{1-1} \hookrightarrow G_1~. \] Let's define them very
briefly. We fix $(X,u)$ in $\sP$ as before and a
$(0,\dots,0)$-simplex $W\to W'\to \sO(X)$ in $G_1(X,u)$.

\begin{defn} To promote $W\to W'\to \sO(X)$ to a $(0,0,\dots,0)$-simplex
in $G_{1-1}(X,u)$, we impose the condition that both $W\to \sO(X)$
and $W'\to \sO(X)$ be simultaneously transverse to the
stratification of $\sO(X)\smin 0$ by strata $\sO(X,\sigma)$, and to
the norm levels. (Compare definition~\ref{defn-G_01}.)
\end{defn}

\begin{defn} To promote $W\to W'\to \sO(X)$ further to a $(0,\dots,0)$-simplex
in $G_{1-2}(X,u)$ we add the following data: contravariant
CW-functors $W_\ddiamond$ and $W'_\ddiamond$ from $\{0\}*\cat(X)$ to
compact spaces, a CW-embedding $W_\ddiamond\to W'_\ddiamond$ and
natural transformations $\gamma\colon W_\ddiamond\to W_\diamond$~,
$\gamma'\colon W'_\ddiamond \to W'_\diamond$ which evaluate to
homotopy equivalences for every object of $\{0\}*\cat(X)$. We
require commutativity of
\[
\xymatrix{
W_\ddiamond \ar[r] \ar[d] & W'_\ddiamond \ar[d] \\
W_\diamond \ar[r] & W'_\diamond\,. }
\]
\end{defn}

\begin{defn}
\label{G1transCWandall} To promote $\colon W\to W' \to \sO(X)$
further to a $(0,0,\dots,0)$-simplex in $G_{1-3}(X,u)$, we add
compatible cellular diagonal approximations for $W_\ddiamond$ and
$W'_\ddiamond$~, and fundamental cycles/chains in the cellular chain
complex(es) of $W_\ddiamond$.
\end{defn}

Similarly as in the case of the functor $G_{0-3}$ the additional
data allow us to extract certain algebraic data. These are two
preferred structures of $n$-dimensional SAPPs (where $n = \dim (W) =
\dim (W')$) and a map
\[ q_\ddiamond \colon (C_*(W_\ddiamond(0)),C_*(W_{\ddiamond}[1]))
\ra (C_*(W'_\ddiamond(0)),C_*(W'_{\ddiamond}[1])) \]
of SAPPs in $\BB(\AA)$, refined to a map of dissected $(n-1)$-dimensional
SAPCs on the boundary in $\BB(\AA_*(X))$. To obtain a single SAPP
with a contractibility property, which is our goal, we need the
construction of \emph{symmetric kernels} in the setting of
chapter~\ref{alg-sur}, definitions~\ref{defn-chduality}
and~\ref{defn-SAPC}. This is a purely algebraic and functorial
construction and is given after Definition
\ref{defn-thegoodtransleft} below.

\emph{Remark.} A $(0,0,\dots,0)$-simplex in $G_{1-i}(X,u)$ is still
considered trivial if the corresponding map $W\to W'$ is a
homeomorphism. All trivial $(0,0,\dots,0)$-simplices are to be
identified with each other.

\begin{lem} The forgetful maps
\[  G_{1-3}\to
G_{1-2} \to G_{1-1} \hookrightarrow G_1
\]
are weak homotopy equivalences.
\end{lem}

All the above refinements of $G_1$ admit actions by the
corresponding refinements of $G_0$ (which refine the action of $G_0$
on $G_1$). At the top of the range we get multiplication and action
maps
\[ G_{0-3}(X,u)
\wedge G_{0-3}(Y,v) \lra G_{0-3}(X*Y,u*v) \]
\[ G_{0-3}(X,u)
\wedge G_{1-3}(Y,v) \lra G_{1-3}(X*Y,u*v).
\]
for any $(X,u)$ and $(Y,v)$ in $\sP$.

\medskip
That completes our efforts to extract chain complex algebra from the
functors $G_0$ and $G_1$. Now we need to do some more work on the
geometric side.

\begin{defn}
\label{defn-norms} Let $s\in [0,1]$. For $V$ in $\sJ\iso$ and $f\co
X\to V$ in $\sP_V$, we define a ``norm'' function on $\sO(X)$ by $tx
\mapsto \|tx\|_s= (1-s)\|tx\|+s\|f(tx)\|$~, using the Euclidean norm
on $V$.
\end{defn}

Keeping the notation of definition~\ref{defn-norms}, we introduce a
refinement $r_fG_0(X,u)$ of $G_0(X,u)$, and a refinement
$r_fG_1(X,u)$ of $G_1(X,u)$, both depending on $f\co X\to V$.

\begin{defn}
\label{defn-rutsch0} Let $V$ be in $\sJ\iso$ and $f\co X\to V$ in
$\sP_V$. To promote a nontrivial $(0,0,\dots,0)$-simplex $p\co W\to \sO(X)$ in
$G_0(X,u)$ to the status of a nontrivial $(0,0,\dots,0)$-simplex in
$r_fG_0(X,u)$ we add the assumption $W=V$ and the following data:
\begin{itemize}
\item a homotopy $(p_s)_{0\le s\le 1}$ through $\ZZ/2$-maps $p_s\co W\to \sO(X)$,
each having preimage of base point equal to base point, such that
$p_0=p$ and $p_s$ is transverse to the nonzero levels of the norm
function $\|...\|_s$ on $\sO(X)$, for $s\in [0,1]$;
\item a homotopy $(h_t)_{1\le t\le 2}$ from the composition
\[
\xymatrix{ W \ar[r]^{p_1} & \sO(X) \ar[r]^f  & V }
\]
to the identity, where each $h_t\co W\to V$ is equivariant, with
preimage of base point equal to base point, and transverse to the
nonzero levels of the euclidian norm on $V$.
\end{itemize}
\end{defn}

\begin{defn}
\label{defn-rutsch1} Let $V$ be in $\sJ\iso$ and $f\co X\to V$ in
$\sP_V$. To promote a $(0,0,\dots,0)$-simplex
\[
\xymatrix{ W \ar[r]^q &  W' \ar[r]^p & \sO(X) }
\]
in $G_1(X,u)$ to the status of a $(0,0,\dots,0)$-simplex in
$r_fG_1(X,u)$ we add the assumption $W'=V$ and the following data:
\begin{itemize}
\item a homotopy $(p_s)_{0\le s\le 1}$ through $\ZZ/2$-maps $p_s\co W'\to \sO(X)$,
each having preimage of base point equal to base point, such that
$p_0=p$ and $p_s$ is transverse to the nonzero levels of the norm
function $\|...\|_s$ on $\sO(X)$, for each $s$;
\item a homotopy $(h_t)_{1\le t\le 2}$ from the composition
\[
\xymatrix{ W' \ar[r]^{p_1} & \sO(X) \ar[r]^f  & V }
\]
to the identity, where each $h_t\co W'\to V$ is equivariant, with
preimage of base point equal to base point, and transverse to the
nonzero levels of the euclidian norm on $V$~;
\item a homotopy $(q_s)_{0\le s\le 2}$ through $\ZZ/2$-maps $W\to W'$~,
each having preimage of base point equal to base point, with
$q_0=q$, such that $p_sq_s$ is transverse to the nonzero levels of
the norm function $\|...\|_s$ on $\sO(X)$ for $s\le 1$, and $h_sq_s$
is transverse to the nonzero levels of the euclidian norm on $V$ for
$s\ge 1$.
\end{itemize}
Such a simplex is \emph{trivial} if $q_s$ is a homeomorphism for all
$s\in[0,2]$. All trivial simplices are to be identified with each other.
\end{defn}

\medskip
Next, $r_fG_{i-3}(X,u)$ is defined (for $i=0,1$) by means of a
pullback square
\[
\xymatrix{ {r_fG_{i-3}(X,u)}
\ar[r]\ar[d] & {G_{i-3}(X,u)} \ar[d]^{\textup{forget}} \\
{r_fG_i(X,u)} \ar[r]^{\textup{forget}} & {G_i(X,u)}. }
\]
To be more precise we form the pullback at the level of
multisimplicial sets, i.e., before realization. The square is also a
homotopy pullback square (by a direct check on homotopy groups). For
$V$ in $\sJ\iso$ let
\[ E^{ga}(V) =
\hocolimsub{f\co X\to V\textup{ in }\sP_V} r_fG_{0-3}(X,u), \]
\[ F^{ga}(V) =
\hocolimsub{f\co X\to V\textup{ in }\sP_V} r_fG_{1-3}(X,u). \] We
still have multiplication and action maps
\[
\begin{array}{ccc}
E^{ga}(V)
\wedge E^{ga}(W) & \lra & E^{ga}(V\oplus W) \\
E^{ga}(V) \wedge F^{ga}(W) & \lra & F^{ga}(V\oplus W)
\end{array}
\]
for $V,W$ in $\sJ\iso$.

\medskip
Let $F^{a,K}$ and $E^{a,K}$ be the functors on $\sJ\iso$ which we
denoted by $F^K$ and $E^K$ in section~\ref{proofs}. We define $F^g$
essentially as the restriction to $\sJ\iso$ of the functor $F$ on
$\sJ$ constructed in \cite{Ma}, except for the small change that
$F^g(W)$ is defined as the geometric realization of an $m$-fold
simplicial set, where $m= \dim(W)$. To make up for the information
lost in restricting from $\sJ$ to $\sJ\iso$, we introduce a
multiplicative functor $E^g$ on $\sJ\iso$ as in example~\ref{expl-terminalaction}.

\begin{defn}
\label{defn-bigEg} For $V$ in $\sJ\iso$, let $E^g(V)=\SS^0$. In each
multidegree $(k_1,\ldots, k_n)$ we have a base-point and another
point which we think of as represented by the space $V \times \prod
\Delta^{k_i}$. The action map $E^g \wedge F^g \ra F^g$ is given by
multiplying a map $q \co W \times \prod \Delta^{l_i} \ra W' \times
\prod \Delta^{l_i}$ with the identity on $V \times \prod
\Delta^{k_i}$.
\end{defn}

\medskip
Now we can relate the functor pair $(E^{a,K},F^{a,K})$ to the pair
$(E^g,F^g)$, using $(E^{ga},F^{ga})$ as a stepping stone. Namely,
there are forgetful natural transformations
\[
\xymatrix@C=40pt@R=2pt{
E^g & \ar[l]_{w_0} E^{ga} \ar[r]^{v_0} & E^{a,K} \\
F^g & \ar[l]_{w_1} F^{ga} \ar[r]^{v_1} & F^{a,K} }
\]
(details in definitions~\ref{defn-thegoodtransleft}
and~\ref{defn-thegoodtransright} below) respecting the
multiplications and the actions. The two in the lower row, $w_1$ and
$v_1$~, are natural homotopy equivalences, giving us (at last) a
two-step identification of $F^g$ with $F^{a,K}$. Of the two in the
upper row, $w_0$ is again a natural homotopy equivalence, whereas
$v_0$ is not. But since $v_0$ is multiplicative, we can use it to
let $E^{ga}$ act on $F^{a,K}$. We then do a derived induction along
$w_0$ to obtain what is essentially an extension of $F^{a,K}$ to
$\sJ\supset\sJ\iso$. See example~\ref{expl-terminalaction} and
lemma~\ref{lem-boring}. This closes the gap between
section~\ref{proofs} and section~\ref{proof}, because the main
results of section~\ref{proofs} can be re-interpreted as results
about $F^{a,K}$ as a functor on all of $\sJ$.

\begin{defn}
\label{defn-thegoodtransleft} The transformation $w_0$ is defined by
taking trivial simplices to the basepoint of $\SS^0$~, and
nontrivial simplices to the non-basepoint. The transformation $w_1$
is induced by forgetful maps
\[ r_fG_{1-3}(X,u) \to r_fG_1(X,u) \to F^g(V). \]
In the notation of definition~\ref{defn-rutsch1}, we proceed by
taking a $(0,0,\dots,0)$-simplex (for example) in $r_fG_1(X,u)$~,
consisting of
\[ \CD W @>q>> W' @>p>> \sO(X) \endCD \]
and homotopies $(p_s)$, $(h_s)$, $(q_s)$ to $q_2\co W\to W'=V$ which
is a $(0,0,\dots,0)$-simplex in $F^g(V)$.
\end{defn}

Before describing $v_0$ and $v_1$, we review the construction of
\emph{symmetric kernels} in the setting of chapter~\ref{alg-sur},
definitions~\ref{defn-chduality} and~\ref{defn-SAPC}. Let $\AA$ be
an additive category with chain duality. Let $f\co C\to D$ be a
morphism in $\BB(\AA)$ and let $\psi$ be an $n$-dimensional
symmetric Poincar\ee structure on $C_*$ such that $f_*\psi$ is also
symmetric Poincar\ee. It is well known that in such a case
$(C,\psi)$ must break up, up to suitable homotopy equivalence, into
a sum
\[  (K_*,\varphi)\oplus (D,f_*\psi) \]
of SAPC's. What we need is a functorial construction of
$(K,\varphi)$. Let $K$ be the mapping cone of the composite chain
map of degree $n$ defined as
\[
\xymatrix@C=40pt{ TD \ar[r]^{Tf} & TC \ar[r]^{\psi_0} & C~. }
\]
(Think of this as an ordinary chain map of degree zero from a
shifted copy of $TD$ to $C$, where the differentials $d\co TD_i\to
TD_{i-1}$ have been multiplied by $(-1)^{ni}$.) Then there is an
inclusion $e\co C\to K$. We let $\varphi:=e_*\psi$~, which is a
symmetric Poincar\ee structure on $K$. The construction works also
for pairs.

Hence, in the notation of comments after Definition
\ref{G1transCWandall}, using the symmetric kernel construction, we
obtain from a map of SAPPs
\[
q_\ddiamond \colon (C_*(W_\ddiamond(0)),C_*(W_{\ddiamond}[1])) \ra
(C_*(W'_\ddiamond(0)), C_*(W'_{\ddiamond}[1]))
\]
an $n$-dimensional SAPP
\[ (K(q_{\ddiamond}(0)),K(q_\ddiamond[1])) \]
in $\BB(\AA)$, refined to a dissected SAPC structure on the boundary
in $\BB(\AA_*(X))$. The condition that the map $q$ is of degree $\pm
1$ assures the contractibility condition.

\begin{defn}
\label{defn-thegoodtransright} The transformation $v_0$ is induced
by forgetful maps
\[ r_fG_{0-3}(X,u) \to G_{0-3}(X,u) \to E^a(X,u) \]
where the second one extracts the chain complex data (including
symmetric structures). The transformation $v_1$ is induced by maps
\[ r_fG_{1-3}(X,u) \to G_{1-3}(X,u) \to F^a(X,u) \]
where the first is forgetful and the second extracts the symmetric
kernels from the available chain complex data (including symmetric
structures).
\end{defn}

The definitions of $v_0$, $v_1$ are set up to respect
the multiplication and action maps.

\medskip
\emph{Remark.} There is a slight complication in the definition
of $v_1$~, due to the fact that the symmetric kernels determined
by multisimplices in $G_{1-3}(X,u)$ which we have defined as trivial
are not completely trivial. (They are contractible but
they are not equal to zero.) It seems best to agree that, wherever a multisimplex
in $G_{1-3}(X,u)$ has trivial (multi)faces, the corresponding subcomplexes
of the symmetric kernel determined by that multisimplex must be
collapsed to zero.

\bigskip
We conclude with an explanation of the remark on infinite loop space
structures at the end of the introduction. Each space $F(X,u)$, in
the notation of section~\ref{proofs}, comes equipped with a
structure of (underlying space of a) $\Gamma$-space in the sense of
Segal \cite{Seg}, determined by the direct sum operation in the
categories $\BB(\AA)$ and $\BB(\AA_*(X))$. This structure is clearly
preserved by the multiplication maps
\[  E(X,u) \wedge F(Y,v) \lra F(X*Y,u*v). \]
(We do not need and we do not use a structure of $\Gamma$-space on
$E(X,u)$ here. Informally, one could say that the adjoint map from
$E(X,u)$ to the space of maps from $F(Y,v)$ to $F(X*Y,u*v)$ factors
canonically through the space of $\Gamma$-maps from $F(Y,v)$ to
$F(X*Y,u*v)$.) It follows that $F^{a,K}$ can be refined to a functor
(first on $\sJ\iso$, then on $\sJ$) with values in the category of
group-like $\Gamma$-spaces. By \cite{Seg}, the ``underlying space''
functor from group-like $\Gamma$-spaces to spaces factors through
the category of infinite loop spaces.


\begin{thebibliography}{99}

\bibitem[Dro]{Dro} \textbf{E Dror-Farjoun}, `Homotopy and homology
of diagrams of spaces', in {\em Algebraic Topology} (Seattle, Wash.,
1985), Springer, Berlin 1987, pp. 93--134.

\bibitem[Hat]{Hat} \textbf{A Hatcher}, `Algebraic Topology', Cambridge
University Press, 2002.

\bibitem[KM]{KM} \textbf{A Kock and I Moerdijk}, {\em Spaces with local equivalence
relations, and their monodromy}, Topology Appl. 72 (1996), 47--78.

\bibitem[Ma]{Ma}\textbf{T Macko}, {\em The block structure
spaces of real projective spaces and orthogonal calculus of
functors}, Trans. Amer. Math. Soc. 359 (2007), 349-383.


\bibitem[Qu]{Qu}\textbf{F Quinn}, `A Geometric Formulation of
Surgery Theory', in {\em Topology of Manifolds} (Proc. Inst., Univ
of Georgia, Athens, GA, 1969), Markham, Chicago, IL, 1970, pp.
500-511.

\bibitem[Ra]{Ra}\textbf{A A Ranicki}, {\em Algebraic L-theory and Topological
Manifolds}, (Cambridge Tracts in Math. 102, Cambridge University
Press, 1992.)

\bibitem[RaLMS1]{RaLMS1}\textbf{A A Ranicki}, {\em The algebraic theory
of surgery, I. Foundations}, Proc. Lond. Math. Soc. 40 (1980),
87--192.

\bibitem[RaLMS2]{RaLMS2}\textbf{A A Ranicki}, {\em The algebraic theory
of surgery, II. Applications to topology}, Proc. Lond. Math. Soc. 40
(1980), 193--287.

\bibitem[RaWe]{RaWe} \textbf{A A Ranicki and M Weiss}, {\em On the algebraic
$L$-theory of $\Delta$-sets}, preprint available from arXiv.

\bibitem[Seg]{Seg} \textbf{G Segal}, {\em Categories and cohomology theories},
Topology 13 (1974), 293-312.

\bibitem[Wa]{Wa}\textbf{C T C Wall}, {\em Surgery on Compact Manifolds},
(2nd edition, edited by \textbf{A A Ranicki}, Mathematical Surveys
and Monographs 69, AMS, 1999.)

\bibitem[We]{We}\textbf{M Weiss}, `Orthogonal Calculus', {\em
Trans. Amer. Math. Soc} 347 (1995), 3743-3796. Erratum, Trans. Amer.
Math. Soc. 350 (1998), 851-855.

\bibitem[We2]{We2}\textbf{M Weiss}, `Visible L-theory', {\em Forum
Mathematicum} 4 (1992), pp. 465-498.



\end{thebibliography}
\end{document}